\documentclass[11pt]{article}%
\usepackage{makeidx}
\usepackage{amssymb}
\usepackage[all]{xy}
\usepackage[latin1]{inputenc}
\usepackage{amsfonts}
\usepackage{amsmath}
\usepackage{amssymb}

\usepackage{url}
\usepackage{hyperref}
\usepackage{graphicx}%
\usepackage{mathabx}

\usepackage{amsthm}
\usepackage{url}
\usepackage{hyperref}
\providecommand{\U}[1]{\protect\rule{.1in}{.1in}}

\usepackage{xcolor}

\newtheorem{prop}{Proposition}[section]
\newtheorem{cor}[prop]{Corollary}

\newtheorem{lem}[prop]{Lemma}

\newtheorem{theo}[prop]{Theorem}

\newcommand{\EE}{\mathbb{E}}

\newcommand{\LL}{\mathbb{L}}

\newcommand{\PP}{\mathbb{P}}
\newcommand{\QQ}{\mathbb{Q}}
\newcommand{\RR}{\mathbb{R}}

\newcommand{\Aa}{ {\cal A }}

\newcommand{\Ca}{ {\cal C }}

\newcommand{\La}{ {\cal L }}
\newcommand{\Na}{ {\cal N }}

\newcommand{\Ea}{ {\cal E }}

\newcommand{\Ra}{ {\cal R }}
\newcommand{\Va}{ {\cal V }}
\newcommand{\Ua}{ {\cal U }}
\newcommand{\Fa}{ {\cal F }}
\newcommand{\Ga}{ {\cal G }}
\newcommand{\Qa}{ {\cal Q }}

\newcommand{\Xa}{ {\cal X }}
\newcommand{\Ma}{ {\cal M }}

\newcommand{\Pa}{ {\cal P }}
\newcommand{\Za}{ {\cal Z }}
\newcommand{\Ya}{ {\cal Y }}

\newcommand{\point}{\mbox{\LARGE .}}

\newcommand{\cqfd}{\hfill\blbx \\}
\def\blbx{\hbox{\vrule height 5pt width 5pt depth 0pt}\medskip}

\def \PP{\mathbb{P}}
\def \RR{\mathbb{R}}
\def \SS{\mathbb{S}}
\def \EE{\mathbb{E}}
\def \QQ{\mathbb{Q}}

\def \LL{\mathbb{L}}

\def \WW{\mathbb{W}}

\def\tr{\mbox{\rm tr}}  

\begin{document}

% =======================================================
  \title{On the stability and the uniform propagation of chaos of a Class of Extended Ensemble Kalman-Bucy filters }
  \author{P. Del Moral, A. Kurtzmann, J. Tugaut}

% =======================================================

\maketitle

\begin{abstract}
This article is concerned with the exponential stability and the uniform propagation of chaos properties of a class of Extended Ensemble
Kalman-Bucy filters with respect to the time horizon. 
This class of nonlinear filters can be interpreted as the conditional
expectations of nonlinear McKean-Vlasov type diffusions with respect to the observation process. We consider filtering problems with 
Langevin type signal processes
 observed by some noisy linear and Gaussian type sensors.
In contrast with more conventional  Langevin nonlinear drift type  processes, the mean field interaction is encapsulated in the covariance matrix of the diffusion. The main results discussed in the article are quantitative estimates of the 
exponential stability properties of these nonlinear diffusions.  These stability properties are used to derive
uniform and non asymptotic estimates of the propagation of chaos properties of Extended Ensemble Kalman filters,
including exponential concentration inequalities. 
To our knowledge these results seem to be the first results of this type for this class of nonlinear
ensemble type Kalman-Bucy filters.\\

\emph{Keywords} : Extended Kalman-Bucy filter, Ensemble Kalman filters, Monte Carlo methods, mean field particle systems,
stochastic Riccati matrix equation, propagation of chaos properties, uniform estimates.\newline

\emph{Mathematics Subject Classification} :  60J60, 60J22, 35Q84, 93E11, 60M20, 60G25.

\end{abstract}

%\tableofcontents
\section{Introduction}

From the probabilistic viewpoint, the Ensemble Kalman filter ({\em abbreviated EnKF}) proposed by G. Evensen in the beginning of the 1990s~\cite{evensen-intro} is
a mean field particle interpretation of extended Kalman type filters. More precisely, Kalman type filters (including the conventional Kalman filter and extended Kalman filters) 
can be interpreted as the conditional expectations of a McKean-Vlasov type nonlinear diffusion. The key idea is to approximate the Riccati equation by a sequence of sample covariance matrices associated with a series of interacting Kalman type filters. 

In the linear Gaussian case these particle type filters converge to the optimal Kalman filter as the number of samples (a.k.a. particles) tends to $\infty$. Little is known for nonlinear and/or non Gaussian filtering problems, apart that they do not converge to the desired optimal filter. This important problem is rather well known in signal processing community. 
For instance, we refer the reader to~\cite{kelly,legland} for a more detailed discussion on these questions in discrete time settings. In this connection, we mention that these ensemble Kalman type filters differ from 
interacting jump type particle filters and related sequential Monte Carlo methodologies. These mean field particle methods
are designed to approximate the conditional distributions of the signal given the observations. It is clearly not the scope of this article to give a comparison between these two different particle methods. For a more thorough discussion we refer the reader to the book~\cite{mean-field} and the references therein.
We also mention that the EnKF models discussed in this article slightly differ from more conventional
EnKF used to approximate nonlinear filtering problems. To be more precise we design a new class of EnKF that converges to
the celebrated extended Kalman filter as the number of  particles goes to $\infty$.
 
These powerful Monte Carlo methodologies are used with success in a variety of scientific disciplines, and more particularly in data assimilation method  for filtering high dimensional problems arising in fluid mechanics and
geophysical sciences~\cite{kalnay,lisa,majda,mandel,nydal,ott,seiler,skj,weng}. A more thorough discussion on the origins and the application domains of EnKF is provided in the series of articles~\cite{burgers,dt-1-2016,evensen-review,evensen-reservoir} and in the seminal research monograph by G. Evensen~\cite{evensen-book}.

The mathematical foundations and the convergence of the EnKF have started in 2011 with the independent pioneering works of 
F. Le Gland, V. Monbet and V.D. Tran~\cite{legland}, and the one by J. Mandel, L. Cobb, J. D. Beezley~\cite{mandel}.
These articles provide $\LL_{\delta}$-mean error estimates for discrete time EnKF and show that they converge towards the Kalman filter as the number of samples tends to infinity. We also quote the recent article by 
D.T. B. Kelly, K.J. Law, A. M. Stuart~\cite{kelly} showing the consistency of Ensemble Kalman filters in continuous and discrete time settings. In the latter the authors show that the Ensemble Kalman filter is well-posed and the mean error variance does not blow up faster than exponentially. The authors also apply a judicious variance inflation technique to strengthen the contraction properties of the Ensemble Kalman filter. We refer to the pioneering article by J.L. Anderson~\cite{anderson-2,anderson,anderson-3} on adaptive covariance inflation techniques, and to the discussion given in the end of Section~\ref{statement-section} in the present article.

In a more recent study by X. T. Tong, A. J. Majda and D. Kelly~\cite{tong}
the authors analyze the long-time behaviour and the ergodicity of discrete generation  EnKF using Foster-Lyapunov techniques ensuring that the filter is asymptotically stable w.r.t. any erroneous initial condition. 
These important properties ensure that the EnKF has a single invariant measure 
and initialization errors of the EnKF will not dissipate w.r.t. the time parameter.
Beside the importance of these properties, the only ergodicity of the particle process does not give any information on the convergence and the accuracy of 
the particle filters towards the optimal filter nor towards any type of extended Kalman filter, as the number of samples tends to infinity.

Besides these recent theoretical advances, the rigorous mathematical analysis of long-time behaviour of these particle methods is still at its
 infancy. As underlined by the authors in~\cite{kelly}, many of the algorithmic innovations associated with the filter, which are required to make a useable algorithm in practice, are derived in an {\em ad hoc} fashion. The divergence of ensemble Kalman filters has been observed numerically in some situations~\cite{georg,kelly-majda,majda-2},
 even for stable signals. This critical phenomenon, often referred as the {\em catastrophic filter divergence} in data assimilation literature, is poorly understood from the mathematical perspective.  
 Our objective is to better understand the long-time behaviour of ensemble Kalman type filters from a mathematical perspective.  
 Our stochastic methodology combines
 spectral analysis of random matrices with recent developments in concentration inequalities, coupling theory and contraction inequalities w.r.t. Wasserstein metrics. 

These developments have been started in two recent articles~\cite{dt-1-2016,dkt-1-2016}. 
The first one provides
uniform propagation of chaos properties of ensemble Kalman filters in the context of linear-Gaussian filtering problems. 
The second article is only concerned with extended Kalman-Bucy filters.
It discusses the stability properties of 
these filters in terms of
 exponential concentration inequalities. These concentration inequalities allow to design confidence intervals
 around the true signal and extended Kalman-Bucy filters. Following these studies, we consider filtering problems with uniformly stable signal processes.

 This condition on the signal is a necessary and sufficient condition to derive uniform estimates for any type of particle filters~\cite{dm-springer2004,mean-field,dt-1-2016} w.r.t. the time parameter. For instance when the sensor matrix is null or for a single particle
 any Ensemble type filter reduces to an independent copy of the signal. In these rather elementary cases, the stability of the signal
 is required to have any type of uniform estimate for any size of the systems.

 We illustrate these models and our stability and observability conditions with a class of nonlinear Langevin type filtering problems, with several classes of sensor models
  
The first contribution of the article is to extend these results as the level of the
McKean-Vlasov type nonlinear diffusion associated with the ensemble Kalman-Bucy filter. Under some natural regularity conditions we show that
these nonlinear diffusions are exponentially stable, in the sense that they forget exponentially fast
any erroneous initial condition. These stability properties are analyzed using coupling techniques and
expressed in terms of 
$\delta$-Wasserstein
metrics. 

The main objective of the article is to analyze the long-time behaviour of the mean field particle interpretation of these nonlinear diffusions. 
We present new uniform estimates w.r.t. the time horizon for the bias and the propagation of chaos properties of the mean field systems. We also quantify
the fluctuations of the sample mean and covariance particle approximations. 

The rest of the article is organized as follows: 

Section~\ref{description-model-sec}
presents the nonlinear filtering problem discussed in the article, the Extended Kalman-Bucy filter,
 the associated nonlinear McKean-Vlasov diffusion and its mean field particle interpretation.
The two main theorems of the article are described in Section~\ref{statement-section}.
In a preliminary short section, Section~\ref{preliminary-section}, we show that the conditional expectations and 
the conditional covariance matrices of  the nonlinear McKean-Vlasov diffusion coincide with the EKF.
We also provide a pivotal fluctuation theorem on the time evolution of these conditional statistics.
Section~\ref{sec-stability} is mainly concerned with the stability properties of the nonlinear 
diffusion associated with the EKF. Section~\ref{pc-section} is dedicated to the propagation of chaos properties
of the extended ensemble Kalman-Bucy filter.

\subsection{Some basic notation}
This section provides with some notation and terminology used in several places in the article.

Given some random variable $Z$ with some probability measure $\mu$ and some function
$f$ on some product space $\RR^r$, we let $$\mu(f)=\EE(f(Z))=\int~f(x)~\mu(dx)$$ be the integral of $f$ w.r.t. $\mu$ or the expectation of $f(Z)$. This notation is rather standard in probability theory. It extends
to integral on Euclidian state spaces the conventional vector summation notation $\mu(f)=\sum_{x} \mu(x)~f(x)$ between row vector measures $\mu=(\mu(x))_{x\in E}$ and 
dual column vector functions $f=(f(x))_{x\in E}$ on finite state spaces $E=\{1,\ldots,d\}$, for some parameter $d\geq 1$.

We let $\left\Vert\point\right\Vert$ be the Euclidean norm on $\RR^{r}$, for some $r\geq 1$. We denote by $\SS_r$ the set of $(r\times r)$ symmetric matrices
with real entries, and by $\SS_r^+$ the subset of positive definite matrices.

We denote by $\lambda_{\tiny min}(S)$ and $\lambda_{\tiny max}(S)$ the minimal and the maximal eigenvalue of  a given symmetric matrix $S$. We let $\rho(P)=\lambda_{\tiny max}((P+P^{\prime}))/2$ be the logarithmic norm of a given square matrix $P$.
Given  $(r_1\times r_2)$ matrices $P,Q$ we define the Frobenius inner product
$$
\langle P,Q\rangle=\tr(P^{\prime}Q)\quad\mbox{\rm and the associated norm}\quad
\Vert P\Vert_F^2=\tr(P^{\prime}P)
$$
where $\tr(C)$ stands for the trace of a given matrix $C$. We also equip the product space $\RR^{r_1}\times \RR^{r_1\times r_1}$ with the inner product 
$$
\langle (x_1,P_1),(x_2,P_2)\rangle:=\langle x_1,x_2\rangle+\langle P_1,P_2\rangle\quad
\mbox{\rm and the norm}\quad \Vert (x,P)\Vert^2:=\langle (x,P),(x,P)\rangle.
$$

Given some $\delta\geq 1$, the $\delta$-Wasserstein distance $\WW_{\delta}$ between two probability measures
$\nu_1$ and $\nu_2$ on some normed space $(E,\Vert\point\Vert)$ is defined by
$$
\WW_{\delta}(\nu_1,\nu_2)=\inf{\EE\left(\Vert Z_1-Z_2\Vert^{\delta}\right)^{1/\delta}
}.
$$
The infimum in the above displayed formula is taken of all pair of random variable $(Z_1,Z_2)$
such that $\mbox{\rm Law}(Z_i)=\nu_i$, with $i=1,2$.

 In the further development of the article, to avoid unnecessary repetitions we also use the letter
$``c"$ to denote some finite constant whose values may vary from line to line, but they do not depend on the time parameter.

\subsection{Description of the models}\label{description-model-sec}

Consider a time homogeneous nonlinear filtering problem of the following form
\begin{equation}\label{filtering-model-ref}
\left\{
\begin{array}{rcl}
dX_t&=&A(X_t)~dt~+~R^{1/2}_{1}~dW_t\\
dY_t&=&BX_t~dt~+~R_2^{1/2}~dV_t
\end{array}
\right.\quad\mbox{\rm and we set $\Ga_t=\sigma\left(Y_s,~s\leq t\right)$. }
\end{equation}
In the above display,  $(W_t,V_t)$ is an $(r_1+r_2)$-dimensional Brownian motion, $X_0$ is a $r_1$-valued
random vector with mean and covariance matrix $(\EE(X_0),P_0)$ (independent of $(W_t,V_t)$), the square root factors 
$R^{1/2}_{1}$ and $R^{1/2}_{2}$ of $R_1$ and $R_2$ are invertible, $B$ is an  $(r_2\times r_1)$-matrix, and $Y_0=0$. The 
drift of the signal 
 is a differentiable vector valued function $
A~:~x\in \RR^{r_1}\mapsto A(x)\in \RR^{r_1}
$ with a Jacobian denoted by $\partial A~:~x\in \RR^{r_1}\mapsto A(x)\in \RR^{(r_1\times r_1)}$.

The Extended Kalman-Bucy filter ({\em abbreviated EKF}) and the associated stochastic Riccati equation are defined by the evolution equations
\begin{eqnarray}\label{EKF}
\left\{
\begin{array}{rcl}
d\widehat{X}_t&=&A(\widehat{X}_t)~dt+P_tB^{\prime}~R^{-1}_{2}~\left[dY_t-B\widehat{X}_t~dt\right]\quad\mbox{\rm with}\quad 
\widehat{X}_0=\EE(X_0),\\
\partial_tP_t&=&\partial A(\widehat{X}_t)P_t+P_t~\partial A(\widehat{X}_t)^{\prime}+R-P_tSP_t\quad\mbox{\rm with}\quad (R,S):=(R_1,B^{\prime}R^{-1}_2B).
\end{array}
\right.
\end{eqnarray}
In the above display, $B^{\prime}$ stands for the transpose of the matrix $B$.

We associate with these filtering models the
conditional nonlinear McKean-Vlasov type diffusion process
\begin{equation}
d\overline{X}_t=\Aa\left(\overline{X}_t,\EE[\overline{X}_t~|~\Ga_t]\right)~dt+R^{1/2}_{1}~d\overline{W}_t+\Pa_{\eta_t}B^{\prime}R^{-1}_{2}~\left[dY_t-\left(B\overline{X}_t~dt+R^{1/2}_{2}~d\overline{V}_{t}\right)\right]\label{EEnKF}
\end{equation}
with the nonlinear drift function
$$
\Aa(x,m):=A\left[m
\right]+
\partial A\left[
m
\right]~(x-m).
$$
In the above display  $(\overline{W}_t,\overline{V}_t,\overline{X}_0)$ stands for independent copies of $(W_t,V_t,X_0)$ (thus independent of
 the signal and the observation path), and  
 $\Pa_{\eta_t}$ stands for the covariance matrix
$$
\Pa_{\eta_t}=\eta_t\left[(e-\eta_t(e))(e-\eta_t(e))^{\prime}\right]
\quad\mbox{\rm with}\quad \eta_t:=\mbox{\rm Law}(\overline{X}_t~|~\Ga_t)\quad\mbox{\rm and}\quad
e(x):=x.
$$
The stochastic process defined in (\ref{EEnKF}) is named the Extended Kalman-Bucy diffusion or simply the EKF-diffusion. 
In Section~\ref{preliminary-section}
(see Proposition~\ref{prop-EKF-diffusion-mean-cov}) we will see that the $\Ga_t$-conditional expectation
of the states $\overline{X}_t$ and their $\Ga_t$-conditional covariance matrices coincide with the EKF filter and the Riccati
equation presented in (\ref{EKF}).

The Ensemble Extended Kalman-Bucy filter ({\em abbreviated En-EKF}) coincides with the mean field particle interpretation
of the nonlinear diffusion process (\ref{EEnKF}). 

To be more precise, let
 $(\overline{W}^i_t,\overline{V}^i_t,\xi^i_0)_{1\leq i\leq N}$  be 
 $N$ independent copies of $(\overline{W}_t,\overline{V}_t,\overline{X}_0)$. 
In this notation, the En-EKF is given by the McKean-Vlasov type interacting diffusion process
\begin{equation}\label{fv1-3}
d\xi^i_t=\Aa(\xi^i_t,m_t)~dt+R^{1/2}_{1}d\overline{W}_t^i+p_tB^{\prime}R^{-1}_{2}\left[dY_t-\left(B\xi^i_t~dt+R^{1/2}_{2}~d\overline{V}^i_{t}\right)\right]
\end{equation}
for any $1\leq i\leq N$,
with  the sample mean and
the rescaled particle covariance matrix defined by
\begin{equation}\label{fv1-3-2}
m_t:=\frac{1}{N}\sum_{1\leq i\leq N}\xi_t^i\quad\mbox{and}\quad
p_t:=\left(1-\frac{1}{N}\right)^{-1}~\Pa_{\eta^{N}_t}=\frac{1}{N-1}\sum_{1\leq i\leq N}\left(\xi^i_t-m_t\right)\left(\xi^i_t-m_t\right)^{\prime}
\end{equation}
with the empirical measures $
\eta^{N}_t:=\frac{1}{N}\sum_{1\leq i\leq N}\delta_{\xi^i_t}$.
We also consider the $N$-particle model $\zeta_t=\left(\zeta^i_t\right)_{1\leq i\leq N}$ defined as $\xi_t=\left(\xi^i_t\right)_{1\leq i\leq N}$
by replacing the sample variance $p_t$ by the true variance $P_t$ (in particular we have $\xi_0=\zeta_0$).

When $B=0$ the En-EKF reduce to $N$ independent copies of the diffusion signal. 
In the same vein, for a single particle the covariance matrix is null so that the En-EKF reduces to a single independent copy of
 the signal. In the case $r_1=1$ we have
 \begin{equation}\label{ref-divergence-1particle}
 \EE\left(\Vert m_t-X_t\Vert^2\right)=2~\mbox{\rm Var}(X_t).
 \end{equation}
 
Last but not least when the number of samples $N<r_1$ is smaller than the dimension of the signal,
 the sample covariance matrix  $p_t$ is the sample mean of $N$ matrices
 of unit rank. Thus, it has at least one null eigenvalue. As a result, in some principal directions the EnKF is only 
 driven  by the signal diffusion. For unstable drift matrices the EnKF experiences divergence as it is not corrected by the innovation 
 process.

 In these rather elementary situations, {\em the stability property of the signal is crucial to
 design some useful uniform estimates w.r.t. the time parameter.} The stability
 of the signal is a also a necessary condition to derive uniform estimates for any type of particle filters~\cite{dm-springer2004,mean-field,dt-1-2016}
 w.r.t. the time parameter.
 
 It should be clear from the above discussion that the stability of the signal is a necessary and sufficient condition to 
design useful uniform estimates w.r.t. the time horizon.

As mentioned in the introduction the En-EKF (\ref{fv1-3}) differs from the more conventional one 
defined as above by replacing $\Aa(\xi^i_t,m_t)$ by the signal drift $A(\xi^i_t)$. In this context the resulting 
sample mean will not converge to the EKF but to the filter defined as in (\ref{EKF}) by replacing
$A(\widehat{X}_t)$ by the conditional expectations $\EE\left(A(X_t)~|~\Ga_t\right)$. The convergence
analysis of this particle model is much more involved than the one discussed in this article. The main difficulty comes from  the dependency on the whole conditional distribution of the signal given the observations. We plan to analyze this class of particle filters in a future study.

\subsection{Regularity conditions}

\subsubsection{Langevin-type signal processes}\label{section-langevin-signal}
 In the further development of the article we assume that the  Jacobian matrix of $A$ satisfies the following regularity conditions:
\begin{equation}\label{regularity-condition-varpiAB}
\left\{
 \begin{array}{rcl}
 -\lambda_{\partial A}&:=&\sup_{x\in \RR^{r_1}} \rho(\partial A(x)+\partial A(x)^{\prime})<0\\
 \\
\Vert \partial A(x)-\partial A(y)\Vert&\leq &\kappa_{\partial A}~\Vert x-y\Vert
\quad\mbox{\rm
 for some $\kappa_{\partial A}<\infty$}
 \end{array}\right.
\end{equation}
where $\rho(P):=\lambda_{\tiny max}(P)$ stands for the maximal eigenvalue of a symmetric matrix $P$. In the above display $\Vert \partial A(x)-\partial A(y)\Vert$ stands for the $\LL_2$-norm of the matrix operator $(\partial A(x)-\partial A(y))$, and $\Vert x-y\Vert$ the Euclidean distance between $x$ and $y$.
A Taylor first order expansion shows that 
\begin{equation}\label{regularity-conditon}
(\ref{regularity-condition-varpiAB})\Longrightarrow \langle x-y,A(x)-A(y)\rangle\leq -\lambda_A~\Vert x-y\Vert^2\quad\mbox{\rm with $\lambda_A\geq \lambda_{\partial A}/2> 0$. }
\end{equation} 

The above rather strong conditions ensure the contraction needed to ensure the stability of the EFK~\cite{dkt-1-2016}. For linear systems $A(x)=Ax$, associated with some matrix $A$, the parameters $\lambda_{A}=\lambda_{\partial A}/2$ coincide with the logarithmic norm of $A$. In this situation we show in~\cite[Section 3.1]{dt-1-2016} that the above condition cannot be relaxed to derive uniform estimates of the Ensemble Kalman-Bucy filter.

The prototype of signals satisfying these conditions are multidimensional diffusions
with  drift functions $(A,\partial A)=(-\partial \Va,-\partial^2 \Va)$ associated with a gradient Lipschitz  strongly convex confining potential $ \Va~:~x\in \RR^{r_1}\mapsto \Va(x)\in [0,\infty[$. The logarithmic norm condition (\ref{regularity-condition-varpiAB}) is met as soon as $\partial^2 \Va\geq v~Id$ with $v=2\vert\lambda_{\partial A}\vert$.   Equivalently the smallest eigenvalue $\lambda_{\tiny min}(\partial^2 \Va(x))$ of the Hessian  is uniformly lower bounded by $v$.
In this case  (\ref{regularity-condition-varpiAB}) is met with $\lambda_{\partial A}=v/2$.  

These conditions are fairly standard in the stability theory of nonlinear diffusions, we refer the reader to the review article~\cite{villani}, and the references therein. Choosing 
 $R_1=\sigma_1^2~Id$ and $A=-\beta\partial  \Va$, for some $\beta,\sigma_1\geq 0$ the signal process $X_t$ resumes to a multidimensional
 Langevin-diffusion
\begin{equation}\label{langevin-example}
dX_t=-\beta~ \partial  \Va(X_t)~dt+\sigma_1~dW_t .
\end{equation}
This process is reversible w.r.t. the invariant distribution. 
Let $\mu$ be a probability distribution
on $\RR^{r_1}$ given by
$$
\mu_{\beta}(dx)=\frac{1}{\Za_{\beta}}~\exp{\left(-\frac{2\beta}{\sigma_1^2} \Va(x)\right)}~dx\quad \mbox{\rm with}\quad \Za_{\beta}=\int~\exp{\left(-\frac{2\beta}{\sigma_1^2} \Va(x)\right)}~dx\in ]0,\infty[.
$$
In the above display $dx$ stands for the Lebesgue measure on $\RR^{r_1}$. 
The Lipschitz-continuity condition of the Hessian $\partial^2 \Va$ introduced in (\ref{regularity-condition-varpiAB}) ensures the continuity of the stochastic Riccati equation (\ref{EKF}) w.r.t. the fluctuations around the random states $\widehat{X}_t$. We illustrate this condition with a
nonlinear example given by the function $$
 \Va(x)=\frac{1}{2}~\langle \Qa_1 x, x\rangle+ \langle q,x\rangle+\frac{1}{3}~\langle \Qa_2 x,x\rangle^{3/2}
$$ with some symmetric positive definite matrices $(\Qa_1,\Qa_2)$ and some given
vector $q\in\RR^{r_1}$. In this case we have
\begin{eqnarray*}
\partial  \Va(x)&=&q+\Qa_1x+\langle \Qa_2 x,x\rangle^{1/2}~\Qa_2 x,\\
\partial^2  \Va(x)&=&\Qa_1+\langle \Qa_2 x,x\rangle^{1/2}~\Qa_2+\langle \Qa_2 x,x\rangle^{-1/2}~\Qa_2 xx^{\prime}\Qa_2 .
\end{eqnarray*}
In this situation we have
\begin{equation}\label{example-non-linear}
\Vert \partial^2  \Va(x)-\partial^2  \Va(y)\Vert\leq 2~\Vert \Qa_2\Vert^{3/2}~\Vert y-x\Vert .
\end{equation}
This shows that conditions (\ref{regularity-condition-varpiAB}) are met with the parameters 
$$(\lambda_{\partial A},\kappa_{\partial A})=\beta~\left(2^{-1}\lambda_{\tiny min}(\Qa_1),~2\lambda_{\tiny max}^{3/2}(\Qa_2)\right) .$$ A proof of (\ref{example-non-linear}) is provided in~\cite[Section 6]{dkt-1-2016}.
More generally these regularity conditions also hold if we replace in (\ref{langevin-example}) the parameter $\sigma_1$ by 
any choice of covariance matrice $R_1$. Also observe that the Langevin diffusion associated with the null form $\Qa=0$ coincides with the conventional linear-Gaussian filtering problem discussed in~\cite{dt-1-2016}. Stochastic gradient-flow diffusions of the form
(\ref{langevin-example}) arise in a variety of application domains. In mathematical finance and mean field game theory~\cite{carmona-fouque,fouque-sun} these Langevin models describe the interacting-collective behaviour of $r_1$-individuals. For instance in the Langevin model 
discussed in~\cite{fouque-sun} the state variables $X_t=\left(X^i_t\right)_{1\leq i\leq r_1}$ represent the log-monetary reserves of $r_1$ banks
lending and borrowing to each other.  The quadratic potential function is given by
$$
\langle \Qa_1 x, x\rangle=\sum_{1\leq i\leq r_1}\left(x_i-\frac{1}{r_1}\sum_{1\leq j\leq r_1}~x_j\right)^2\Rightarrow~\Qa_1
\succ \left(1-\frac{1}{r_1}\right)~I_{r_1} .
$$
In this context, the parameter $\beta$ represents the mean-reversion rate between banks. More general interacting potential functions can be considered.  Mean field type diffusion processes are also used to design low-representation of fluid flow velocity fields. These 
vortex-type particle filtering problems are developed in some details in the pionnering articles by E. M\'emin and his co-authors~\cite{cuzol1,cuzol2,cuzol3,papadakis}. These probabilistic interpretations of the 2d-incompressible Navier-Stokes equation 
represent the vorticity map as a mixture of basis functions centered around
each vortex.  

In this connexion, we mention that our approach also applies to interacting diffusion  gradient flows described 
by a potential function of the form
$$
\Va(x)=\sum_{1\leq i\leq r_1}~\Ua_1(x_i)+\sum_{1\leq i\not=j\leq r_1}~\Ua_2(x_i,x_j)
$$
for some gradient Lipschitz  strongly convex confining potential $ \Ua_i:\RR^i\mapsto [0,\infty[$, $i=1,2$. In this situation, we have
\begin{equation}\label{U1U2U}
\partial^2 \Ua_1\geq u_1\quad\mbox{\rm and}
\quad \partial^2\Ua_2\geq u_2~I_2\Longrightarrow~\partial^2\Va\succ v~I_{r_{1}}\quad \mbox{\rm with}\quad v:=(u_1+(r_1-1)u_2)>0.
\end{equation}
We further assume that
\begin{eqnarray*}
\vert \partial^2\Ua_1(x_1)-\partial^2\Ua_1(y_1)\vert&\leq& \kappa_{\partial^2\Ua_1}~\vert x_1-x_2\vert ,\\
\Vert \partial^2\Ua_2(x_1,x_2)-\partial^2\Ua_2(y_1,y_2)\Vert &\leq &\kappa_{\partial^2\Ua_2}~\Vert (x_1,x_2)-(y_1,y_2)\Vert .
\end{eqnarray*}
In this case,  we have
\begin{equation}\label{2nd-formula-Lipschitz}
\Vert \partial^2\Va(x)- \partial^2\Va(y)\Vert\leq \kappa_{\partial^2\Va}~\Vert x-y\Vert\quad\mbox{\rm with}\quad
 \kappa_{\partial^2\Va}:=\kappa_{\partial^2\Ua_1}+\kappa_{\partial^2\Ua_2}~(r_1-1)~ \sqrt{2(r_1-1)} .
\end{equation}
This shows that conditions (\ref{regularity-condition-varpiAB}) are met with $$(\lambda_{\partial A},\kappa_{\partial A})=\beta~\left(2^{-1}(u_1+(r_1-1)u_2),~\kappa_{\partial^2\Ua_1}+\kappa_{\partial^2\Ua_2}~(r_1-1)~ \sqrt{2(r_1-1)}\right) .$$ 
The detailed proofs of (\ref{U1U2U})-(\ref{2nd-formula-Lipschitz}) are provided in~\cite[Section 6]{dkt-1-2016}.

\subsubsection{Observability conditions}
To introduce our observability conditions we give a brief introduction to the class of observation processes discussed in this article.
When the observation variables are the same as the ones of the signal, the signal observation has the same dimension as the signal and resumes to
 some equation of the form
\begin{equation}\label{fully-observed-sensor}
dY_t=b~ X_t~dt+\sigma_2~dV_t
\end{equation}
for some parameters $b\in\RR$ and $\sigma_2\geq 0$.  These sensors are used in data grid-type assimilation
problems when measurements can be evaluated at each cell. These fully observed models are discussed  in~\cite[Section 4]{harlim-hunt} in the context of the Lorentz-96 filtering problems. These observation processes are also used in the article~\cite{berry-harlim} for application to nonlinear and multi-scale filtering problem. In this context, the observed variables represent the slow components of the signal.  When the fast components are represented by some Brownian motion with a prescribed covariance matrix, the filtering of the slow components with full observations take the form (\ref{fully-observed-sensor}). 

 For partially observed signals we cannot expect any stability properties of the EKF and the En-EKF without introducing some structural conditions of observability and controllability on the signal-observation equation (\ref{filtering-model-ref}). To get one step further in our discussion, observe that the EKF equation (\ref{EKF}) implies that
 \begin{equation}\label{ref-inro-XEKF}
d(\widehat{X}_t-X_t)=\left[(A(\widehat{X}_t)-A(X_t))-P_tS(\widehat{X}_t-X_t)\right]~dt+P_{t}~C^{\prime}R^{-1/2}_{2}~dV_t+R^{1/2}_{1}~dW_t.
 \end{equation}
 
 This equation shows that the stability properties of this process depends on the nature of the real eigenvalues of the symmetric 
 matrices $(A(x)-P_tS)_{\tiny sym}$, with $x\in\RR^{r_1}$. In contrast with the conventional Kalman-Bucy filter, the Riccati equation 
 (\ref{EKF})  is a stochastic equation. 
 
 In this connection, we already mentioned that the sample covariance matrices $p_t$ of the En-EKF also satisfy
 a stochastic Riccati type equation of the same form, up to some fluctuation martingale (see for instance (\ref{f23}) in Theorem~\ref{fluctuation-theorem-ref} in the present article). In the same vein, we shall see in (\ref{f21}) that the En-EKF sample mean $m_t$ evolution satisfies 
the same equation as the EKF, up to some fluctuation martingales coming from the fluctuations of the sample-covariance matrices  and the ones of the sample-particles. 
 
 As a result, the stability properties of the EKF and the En-EKF are not
induced by some kind of observability condition that ensures the existence of a steady state deterministic covariance matrix. The random fluctuations of the matrices $\partial A(\widehat{X}_t)$ and $\partial A(m_t)$  as well as the fluctuations of the stochastic matrices $(A(m_t)-p_tS)_{\tiny sym}$ may corrupt the stability in the EKF and the En-EKF, even if the linearized filtering problem around some chosen state is observable and controllable.  For instance the empirical covariance matrices may not be invertible for small sample sizes. For a more thorough discussion on the stability properties of Kalman-Bucy filters, the EKF and Riccati equations we refer the reader to~\cite{dt-1-2016,dkt-1-2016}, and the references therein.

 As shown in the system above, 
these  fluctuations enter in two different ways in the En-EKF. The first one in the drift function of the system, the other one through the diffusive part.

Therefore the fluctuations of the empirical
 covariances from small sample sizes  corrupt the natural stabilizing effect of the observation process in the EKF filter evolution. 
   In practice it has been observed that these fluctuations induce an underestimation of the true error covariances. As a result the En-EKF eventually
   ignores the information given by the observations. This lack of observation-driven component also leads to the divergence of the filter.
 
Last but not least, from another numerical viewpoint, the En-EKF is also known to be not robust, in the sense that arithmetic errors may accumulate even if the exact filter is stable.
 
 All of these instability properties of the EnKF are well-known and often referred as the catastrophic filter divergence in data assimilation literature, see for instance~\cite{georg,kelly-majda,majda-2}, and the references therein. As mentioned by the authors in~\cite{kelly-majda}, ``catastrophic filter divergence is a well-documented but mechanistically mysterious phenomenon whereby ensemble-state estimates explode to machine infinity despite the true state remaining in a bounded region". In all the situations discussed above the instability properties of Ensemble Kalman-Bucy type filters
 are related to some observability problem. 
 
 The stability analysis of diffusion processes is always much more documented than the ones on their possible divergence.
For instance, in contrast with conventional Kalman-Bucy filters, the stability properties of the EnKF are not
induced by some kind of observability or controllability condition. The only known results for discrete generation
EnKF is the recent work by X. T. Tong, A. J. Majda and D. Kelly~\cite{tong}. One of the main assumptions of the article is that the sensor-matrix 
has full rank. The authors also provide a concrete numerical example of filtering problem with sparse observations for which the EnKF experiences a catastrophic divergence. These divergence properties have been analysed in some details in the article~\cite{dt-1-2016} in the context of linear-Gaussian filtering problems. The full rank observation assumption avoids the EnKF to experience local or global exponential instabilities.

 To quantify and control uniformly in time the propagations of these instabilities we need to introduce
 some strong observability condition that ensures that the system is globally and locally stable.
In the further development of the article we assume that the following condition is satisfied:
\begin{equation}\label{new-condition}
\hskip-5cm(\mbox{\rm S})\hskip3cm \quad S=\rho(S)~Id\quad\mbox{\rm for some}\quad \rho(S)>0.
\end{equation}

The fully observed model discussed in (\ref{fully-observed-sensor}) clearly satisfies
condition (\ref{new-condition})  with the parameter $\rho(S)=(b/\sigma_2)^2$. Condition (\ref{new-condition})
ensures that the particle EnKF has uniformly bounded  $\LL_n$-moments for any $n\geq 1$.
In the context of linear-Gaussian filtering problems, 
this condition is also essential to ensure the uniform convergence of Ensemble Kalman-Bucy filter w.r.t. the time parameter~\cite{dt-1-2016}.
This article also provides a geometric description of the divergence regions in the set of positive covariance matrices for elementary $2$-dimensional observable and controllable systems. When condition (S) is not met, we design stochastic observers driven by these matrices that diverge 
when the signal drift matrix is unstable    (see~\cite[Section 4]{dt-1-2016}). 

From the pure mathematical viewpoint the observability condition (S)  allows to 
combine exponential semigroup techniques with spectral analysis and log-norm inequalities. To get some intuition and to better connect this work with~\cite{dt-1-2016} we give some brief comments on these spectral techniques:

For $2$-dimensional linear signals $A(x)=A x$, the existence and the uniqueness of the steady state $P$ of the Riccati equation (\ref{EKF}) is ensured by some appropriate observability and controllability conditions. In this context we have
$\mu(A-PS)<0$ even for unstable signal-drift matrices. 
This condition ensures the stability of the steady state filter.

 Starting from the steady state $P_0=P$ the EnKF filter is driven by stochastic matrices $p_t$ that converge to $P$, as the size $N$ of the ensemble tends to $\infty$. 
The stability analysis of the EnKF filter now depends on the 
sign of the log-norms $t\mapsto\mu(A-p_tS)$ of the stochastic matrices. The fluctuations of $p_t$ around $P$ are defined by the matrices
\begin{equation}\label{true-perturbations}
Q_t:=\sqrt{N}~(p_t-P)\in\SS_{r_1}\Longleftrightarrow p_t=P+\frac{1}{\sqrt{N}}~Q_t\in\SS_{r_1}^+ .
\end{equation}
Under condition (S) we have $\mu(A)<0\Rightarrow\mu(A-p_tS)=\mu((A-PS)-Q_tS)<\mu(A)<0$ for any fluctuation matrices $Q_t$.
When (S) is not met the local divergence domain of matrices
$Q_t$ such that $\mu(A-p_tS)=\mu((A-PS)-Q_tS)>0$ may be very large, even when $\mu(A)<0$. The refined analysis on the stability of these models requires to
analyze in some details the random excursion of the matrices into these local divergence domains. For a more thorough discussion on these
local and global divergence issues in the context of linear systems we refer the reader to~\cite[Section 4]{dt-1-2016}.

Last but not least, we mention that (\ref{new-condition}) is satisfied when the filtering problem is similar to the ones discussed above; that is, up to a change of basis functions. More precisely, any filtering problem (\ref{filtering-model-ref}) with $r_1=r_2$ and s.t.  $(R_2^{-1/2}B)$ is invertible can be turned  into a filtering problem equipped with an identity sensor matrix; even when the original matrix $S=C^{\prime}R_2^{-1}C=\Ca^{\prime}\Ca$ does not satisfy (\ref{new-condition}). 
To check this claim we observe that
$$
\Ya_t:=R^{-1/2}_2Y_t\quad\mbox{\rm and}\quad \Xa_t:=R^{-1/2}_2BX_t
\Longrightarrow \left\{\begin{array}{rcl}
d\Xa_t&=&\Aa(\Xa_t)~dt+\Ra_1^{1/2}~dW_t\\
d\Ya_t&=&\Xa_t~dt+dV_t
\end{array}\right.
$$
with the drift function
$$
\Aa:=(R^{-1/2}_2B)\circ A\circ (R^{-1/2}_2B)^{-1}\quad \mbox{\rm and the matrix}\quad
\Ra_1:=R^{-1/2}_2 BR_1B^{\prime}R^{-1/2}_2.
$$
In this situation the filtering model $(\Xa_t,\Ya_t)$ satisfies (\ref{new-condition}). In addition, we have
$$
A=(R^{-1/2}_2B)^{-1} \circ\partial U \circ (R^{-1/2}_2B)\Rightarrow
(\Aa,\partial \Aa)=(\partial U,\partial ^2U)
$$
%In this situation the filtering model $(\Xa_t,\Ya_t)$ satisfies (\ref{new-condition}) 
and the signal process $\Xa_t$ belongs to the class of Langevin type diffusion  discussed in Section~\ref{section-langevin-signal}.

\section{Statement of the main results}\label{statement-section}
\subsection{Concentration inequalities}
One of our results concerns the stability properties of the EKF-diffusion (\ref{EEnKF}). 
It is no surprise that these properties strongly 
depend on logarithmic norm of the drift function $A$ as well as on the size of covariance matrices of the signal-observation diffusion. For instance, we have the uniform moment estimate
\begin{equation}\label{uniform-estimates-ref}
\lambda_{\partial A}>0\Rightarrow\forall \delta\geq 1\quad
\sup_{t\geq 0}{\left\{\EE[\Vert X_t\Vert^{\delta}]\vee\tr(P_t)\vee\EE[\Vert X_t-\widehat{X}_t\Vert^{\delta}]\right\}}\leq c.
\end{equation}
A detailed proof of these stochastic stability properties including exponential concentration inequalities can be found in~\cite{dkt-1-2016}. Observe that $\tr(P_t)$ is random so that the above inequality provides an almost sure estimate. To be more precise we use (\ref{EKF}) to check that
\begin{equation}\label{riccati-bound}
\partial_t\tr\left(P_t\right)
\leq  -\lambda_{\partial A}~\tr\left(P_t\right)+\tr(R)\Longrightarrow
\tr(P_t)\leq e^{-\lambda_{\partial A}t}~\tr\left(P_0\right)+{\tr(R)}/{\lambda_{\partial A}}.
\end{equation}
The detailed proof of (\ref{riccati-bound}) can be found on page~\pageref{Lemma-Unif-control-sec}.

To get one step further in our discussion, we consider the following ratio
$$
\lambda_{S}:=\frac{\lambda_{\partial A}}{\rho(S)}, \qquad
\lambda_{R}:=\frac{\lambda_{\partial A}}{\tr(R)}\quad\mbox{\rm and}\quad
\lambda_{K}:=\frac{\lambda_{\partial A}}{\kappa_{\partial A}}.
$$
%The quantity $\rho(S)$ is connected to the sensor matrix $B$ and to the inverse of the covariance matrix of the 
%observation perturbations. We also have the rather crude estimate
%$$
%\rho(S)\leq \tr(S)=\Vert\partial R^{-1/2}_2 B\Vert_F^2\leq  \tr(R_2^{-1})~\Vert B\Vert_F^2
%$$
%
 Roughly speaking, the three quantities presented above measure the relative stability index of the signal drift with respect to
the perturbation degree of the sensor, the one of the signal, and the modulus of continuity of the 
Jacobian entering into the Riccati equation. For instance, $\lambda_S$ is high for sensors with
large perturbations, inversely $\lambda_R$ is large for signals with small perturbations. Most of our analysis  relies on the behaviour of the following quantities:
\begin{eqnarray*}
\lambda_{R,S}&:=&(8e)^{-1}~\lambda_{R}~\sqrt{\lambda_{S}}~
  \left[1+\frac{2}{\lambda_{R}\lambda_S}\right]^{-1},
\\
\widehat{\lambda}_{\partial A}/\lambda_{\partial A}&:=&\left(\frac{1}{2}-\frac{2}{\lambda_{K}\lambda_{R}}\right)+\left(\frac{1}{2}-~
\frac{1}{\sqrt{\lambda_{S}}}\right)\left[1-\frac{3}{4}
\frac{1}{\sqrt{\lambda_{S}}}\right]. \end{eqnarray*}

 In the quadratic Langevin-signal 
filtering problem discussed in (\ref{langevin-example}) and (\ref{fully-observed-sensor}) with $b=1=\beta$  these parameters resume to
\begin{equation}\label{ref-examples-lambda}
\lambda_{S}:=\frac{1}{2}~\sigma_2^2~v , \qquad
\lambda_{R}:=\frac{1}{2r_1}~\sigma_1^{-2}~v\quad\mbox{\rm and}\quad
\lambda_{K}:=\infty .
\end{equation}
In this situation we have
\begin{eqnarray*}
\widehat{\lambda}_{\partial A}/\lambda_{\partial A}&:=&\frac{1}{2}+\frac{1}{2}~\left(1-2\sqrt{2}
\frac{1}{\sigma_2~\sqrt{v}}\right)\left[1-\frac{3\sqrt{2}}{4}
\frac{1}{\sigma_2~\sqrt{v}}\right] . \end{eqnarray*}
Notice that these parameters do not depend on the dimension of the signal, nor on the diffusion parameter $\sigma_1$. 
In addition, we have $\widehat{\lambda}_{\partial A}/\lambda_{\partial A}>0$
for any choice of parameters $(v,\sigma_2)$.

To better connect these quantities with the stochastic stability of the EKF diffusion we discuss some exponential concentration inequalities that can be easily derived from our analysis.
These concentration inequalities are of course more accurate than any type of mean square error estimate. 
Let 
$\widehat{X}_t(m,p)$ be the solution of the EKF equation (\ref{EKF}) starting at $(\widehat{X}_0,P_0)=(m,p)$, and let 
$X_t(x)$ be the state of the signal starting at $X_0(x)=x$. Let $\varpi(\delta)$ be the function
$$
\delta\in [0,\infty[\mapsto \varpi(\delta):= \frac{e^2}{\sqrt{2}}~\left[\frac{1}{2}+\left(\delta+\sqrt{\delta}\right)\right].
$$
In this notation, we have the following exponential concentration inequalities.
\begin{theo}
For  any time horizon $t\in [0,\infty[$, and
any $\delta\geq 0$
 the probabilities  of the following
events 
\begin{eqnarray*}
\Vert X_t(x)-\widehat{X}_t(m,p)\Vert^2
\displaystyle&\leq& \frac{1}{2e}~\varpi(\delta)~{\sqrt{\lambda_{S}}}/{\lambda_{R,S}}
\\
&&+2~
e^{-\lambda_{\partial A}t}~\Vert x-m\Vert^2+8~\varpi(\delta)~
\frac{\vert e^{-\lambda_At}- e^{-\lambda_{\partial A}t}\vert}{\vert \lambda_A/\lambda_{\partial A}-1\vert}~
\tr(p)^2 /\lambda_{S}
\end{eqnarray*}
and
\begin{eqnarray*}
\Vert \overline{X}_t(m,p)-\widehat{X}_t(m,p)\Vert^2&\leq &\frac{1}{2e}~\varpi(\delta)~{\sqrt{\lambda_{S}}}/{\lambda_{R,S}}+8~\varpi(\delta)~
e^{-\lambda_{\partial A}t}~\tr(p)^2/\lambda_{S}
\end{eqnarray*}
are greater than $1-e^{-\delta}$.
\end{theo}
The proof of the first assertion is a consequence of~\cite[Theorem 1.1]{dkt-1-2016},
the proof of the second one is a consequence of the $\LL_{\delta}$-mean error estimate (\ref{unif-control}). 
These concentration inequalities show that the quantity
$$
{\sqrt{\lambda_{S}}}/{\lambda_{R,S}}=8e~\lambda_{S}~\frac{1}{\lambda_{R}\lambda_{S}}~\left[1+\frac{2}{\lambda_{R}\lambda_S}\right]
$$
can be interpreted as the size of a confidence interval around the values of {\em the true signal}, as soon as the time horizon is large. 
It is also notable that the same quantity controls the fluctuations of the EKF diffusion around the values
of the EKF. These confidence intervals are small for stable signals with small perturbations. In the quadratic Langevin-signal 
filtering problem discussed in (\ref{langevin-example}) and (\ref{fully-observed-sensor}) with $b=1=\beta$  the above quantity resumes to
$$
{\sqrt{\lambda_{S}}}/{\lambda_{R,S}}=2^4e~\frac{1}{v}~r_1~\sigma_1^2~\left[1+\frac{8}{v^2}~r_1~\left(\frac{\sigma_1}{\sigma_2}\right)^2\right].
$$
For unit signal-to-noise ratio $\sigma_1=\sigma_2$ these fluctuation parameters are small for stable signals with small perturbations. The above formula also indicate the degradation of the fluctuation parameter when the size of the system is large.

\subsection{A stability theorem}
We further assume that
 \begin{equation}\label{reference-stability}
 (\lambda_{K}\lambda_{R}/4)\wedge \lambda_{R,S}\wedge (\lambda_{S}/4)
 >1.
  \end{equation}
This regularity property is a purely technical condition. The condition  $(\lambda_{K}\lambda_{R}/4)\wedge (\lambda_{S}/4)$ 
 ensures that
$
0< \widehat{\lambda}_{\partial A}\leq \lambda_{\partial A}
$, while $\lambda_{R,S}>1$ is used to derive $\LL_p$-mean error estimates with some parameter $p\geq 1$  that depends on $\lambda_{R,S}$.

The condition (\ref{reference-stability}) is clearly met as soon as $\lambda_{R}$ and $\lambda_{S}$ are sufficiently large.
As we shall see the quantity $\widehat{\lambda}_{\partial A}$ represents the Lyapunov stability exponent of the EKF.
This exponent is decomposed into two parts. The first one represents the relative contribution of the signal perturbations, the second one is related to the sensor perturbations.

In contrast with the linear-Gaussian case discussed in~\cite{dt-1-2016}, the stochastic Riccati equation~\eqref{EKF} depends 
on the states of the EKF.
As shown in~\cite{dkt-1-2016} the stability of the EKF
relies on a stochastic Lyapunov exponent that depends on the 
random trajectories of the filter as well as on the signal-observation processes.
The technical condition (\ref{reference-stability}) allows to control uniformly the fluctuations of these stochastic
exponents with respect to the time horizon.

A more detailed discussion on the regularity condition (\ref{reference-stability}), including a series of sufficient conditions are provided in the appendix, Section~\ref{reg-conditions-appendix}. For filtering problems with an observation process of the form (\ref{fully-observed-sensor}) with $\rho(S)=(b/\sigma_2)^2=1$ we have
$$\lambda_S=\lambda_{\partial A}\quad\Longrightarrow\quad
\lambda_{R,S}:=\frac{1}{8e\tr(R)}~
 \frac{\lambda_{\partial A}^{3+1/2}}{\lambda_{\partial A}^2+2\tr(R)}.
$$
In this situation (\ref{reference-stability}) is met as soon as the following easy to check condition is satisfied 
\begin{equation}\label{CS-easy-check}
\lambda_{\partial A}>4\quad\mbox{\rm and}\quad
\tr(R)\leq \frac{\lambda_{\partial A}^2}{2}\left\{\frac{1}{2\kappa_{\partial A}}\wedge\left[\sqrt{1+\frac{1}{4e\sqrt{\lambda_{\partial A}}}}-1\right]\right\} .
\end{equation}
A detailed proof of this assertion is provided in the end of Section~\ref{reg-conditions-appendix}. In the quadratic Langevin-signal 
filtering problem discussed in (\ref{langevin-example}) and (\ref{fully-observed-sensor}) with $b=\sigma_2$, condition (\ref{CS-easy-check}) resumes to
$$
v/8>1\quad\mbox{\rm and}\quad
2\sqrt{2}e~r_1\sigma_1^2\leq (v/8)~\frac{1}{\sqrt{1+\frac{1}{2\sqrt{2}ev}}+1} .
$$
These conditions are clearly much stronger than the ones discussed in~\cite{dt-1-2016} in the context of linear-Gaussian filtering problems. For the same type of filtering problem, exponential stability and uniform propagations of chaos for the EnKF hold as soon as $v>0$.

Let  $(\overline{X}_t,\overline{Z}_t)$ be a couple of EKF Diffusions (\ref{EEnKF}) starting from two random states
with mean $(\widehat{X}_0,\widecheck{X}_0)$ and covariances matrices $(P_0,\widecheck{P}_0)$ (and driven by the same Brownian motions $(\overline{W}_t,\overline{V}_t)$). One key feature of these nonlinear diffusions is that the $\Ga_t$-conditional expectations $(\widehat{X}_t,\widecheck{X}_t)$  and the 
$\Ga_t$-conditional  covariance matrices $(P_t,\widecheck{P}_t)$ satisfy the EKF and the stochastic Ricatti equations discussed in (\ref{EKF}).

Whenever condition (\ref{reference-stability}) is satisfied we recall from~\cite{dkt-1-2016} that
for any $\epsilon\in ]0,1]$  there exists some time horizon $s$ such that
for any $t\geq s$ we have the almost sure contraction estimate 
$$
\EE\left(\Vert(\widehat{X}_t,P_t)-
(\widecheck{X}_t,\widecheck{P}_t)\Vert^{\delta_S}~|~\Ga_s\right)^{2/\delta_S}
\displaystyle\leq  \Za_{s}~\displaystyle \exp{\left[-\left(1-\epsilon\right)\widehat{\lambda}_{\partial A}(t-s)\right]}~\Vert(\widehat{X}_s,P_s)-
(\widecheck{X}_s,\widecheck{P}_s)\Vert^2
$$
with $\delta_S:= 2^{-1}~\sqrt{\lambda_{S}}$, and some random process $\Za_{t}$ satisfying the uniform moment condition
\begin{equation}\label{unif-moments-Za}
\sup_{t\geq 0}{\EE\left( 
\Za^{\alpha}_t
\right)}<\infty
\quad\mbox{with}\quad \alpha=2\lambda_{R,S}~\delta_S.
\end{equation}
These conditional contraction estimates can be used to quantify the stability properties of the EKF. More precisely,
if we set
$$
\PP_t=\mbox{\rm Law}(\widehat{X}_t,P_t)\quad\mbox{\rm and}\quad
\widecheck{\PP}_t=\mbox{\rm Law}(\widecheck{X}_t,\widecheck{P}_t)
$$
then the above contraction inequality combined with the uniform estimates (\ref{uniform-estimates-ref}) readily implies that
$$
\forall t\geq t_0\qquad
\WW_{\delta_S}^2(\PP_t,\widecheck{\PP}_t)\leq c~
\exp{\left[-t~(1-\epsilon)~\widehat{\lambda}_{\partial A}\right]}
$$
 for any $\epsilon\in [0,1[$, with some time horizon $t_0$. This stability property ensures
 that the EKF forgets exponentially fast any erroneous initial condition. Of course these forgetting properties
 of the EKF do not give any information at the level of the process.
 One of the main objective of the article is to
 complement these conditional expectation stability properties at the level of the McKean-Vlasov type
 nonlinear EKF-diffusion (\ref{EEnKF}). 
 
 Our second main result can basically be stated as
 follows.
\begin{theo}\label{th-intro-1}
 Let
$(\overline{\eta}_t,\breve{\eta}_t)$ be the probability distributions of a couple $(\overline{X}_t,\overline{Z}_t)$ of EKF Diffusions (\ref{EEnKF}) starting from two possibly different 
random states.
Assume condition (\ref{reference-stability}) is met with $\delta^{\prime}_S:= \delta_S/4~\geq 2$. In this situation, for any $\epsilon\in [0,1[$ there exists some time horizon $t_0$ such that for any $t\geq t_0$ we have
\begin{equation}\label{w}
\WW_{\delta_S^{\prime}}^2(\overline{\eta}_t,\breve{\eta}_t)\leq c~
\exp{\left[-t~(1-\epsilon)~\lambda\right]}
\quad\mbox{with}\quad
\lambda\geq \widehat{\lambda}_{\partial A}\wedge({\lambda_{\partial A}}/{4}).
\end{equation} 
\end{theo}

\subsection{A uniform propagation of chaos theorem}

Our next objective is to analyze the long-time behaviour of the mean field type En-EKF model 
discussed in (\ref{fv1-3}). From the practical estimation point of view, only the sample mean and the sample covariance
matrices (\ref{fv1-3-2}) are of interest since these quantities converge to the EKF and the Riccati equations, as $N$ tends to $\infty$.  Another important problem
is to quantify the bias of the mean field particle approximation scheme. These properties are related to the propagation of chaos
properties of the mean field particle model. They are expressed in terms of 
 the collection of probability distributions
$$
\PP_t^N=\mbox{\rm Law}(m_t,p_t),\qquad
\QQ_t^{N}=\mbox{\rm Law}(\xi^1_t)\quad\mbox{\rm and}\quad
\QQ_t=\mbox{\rm Law}(\zeta^1_t).
$$

\begin{theo}
Assume that  (\ref{reference-stability}) is met with  $\delta_{R,S}:=  (e\lambda_{R,S})\wedge \delta_S\geq 2$. 
In this situation,
there exist some $N_0\geq 1$ and some $\beta\in ]0,1/2]$ such that for any $N\geq N_0$, 
 we have
 the uniform non asymptotic estimates
\begin{equation}\label{pc}
 \tr(P_0)^2\leq\frac{\lambda_S}{\lambda_R}
 \left[\frac{1}{2}+\frac{1}{\lambda_{R}\lambda_S}\right]\quad\Longrightarrow\quad
\sup_{t\geq 0}{\WW_{\delta_{R,S}}\left(\PP_t^N,\PP_t\right)}\leq c N^{-\beta}.
\end{equation} 
In addition, when  $\delta_{R,S}\geq 4$ we have the uniform propagation of chaos estimate
\begin{equation}\label{pc-2}
 \sup_{t\geq 0}{\WW_{2}\left(\QQ_t^{N},\QQ_t\right)}\leq cN^{-\beta}.
\end{equation} 
\end{theo}

Our analysis does not provide an explicit formula for the rate of convergence $\beta$. We conjecture that the optimal rate is $\beta=1/2$
as in the linear-Gaussian case developed in~\cite{dt-1-2016}. 

For the quadratic Langevin-signal 
filtering model discussed  in (\ref{langevin-example}) and (\ref{fully-observed-sensor}) with $b=1=\beta$, by 
(\ref{ref-examples-lambda})  the l.h.s. condition in (\ref{pc}) resumes to
$$
 \tr(P_0)^2\leq\frac{\lambda_S}{\lambda_R}
 \left[\frac{1}{2}+\frac{1}{\lambda_{R}\lambda_S}\right]=r_1~(\sigma_1\sigma_2)^2
 \left[\frac{1}{2}+r_1~\left(\frac{2}{v}\right)^2~\left(\frac{\sigma_1}{\sigma_2}\right)^2\right].
$$

We end this section with some comments on our regularity conditions.

The condition $\eqref{new-condition}$ is needed to control the fluctuations 
of the trace of the sample covariance matrices of the En-EKF, even if the trace expectation is uniformly stable. We believe that this technical observability condition can be relaxed.

Despite our efforts, our regularity conditions are stronger than the ones
discussed in~\cite{dt-1-2016} in the context of linear-Gaussian filtering problems.  The main difference here is that the signal stability is required to compensate the possible instabilities created by highly informative sensors when 
we initialize the filter with wrong conditions.  

Next we comment the trace condition in (\ref{pc}).
As we mentioned earlier, the stability properties of the limiting EKF-diffusion (\ref{EEnKF}) 
are expressed in terms of a stochastic Lyapunov exponent that depends on the trajectories
of the signal process. The propagation of chaos properties of the mean field particle approximation  (\ref{fv1-3})
depend on the long-time behaviour of these stochastic Lyapunov exponents. Our analysis
is based on a refined analysis of Laplace transformations associated with quadratic type
stochastic exponents. The existence of these  $\chi$-square type Laplace transforms requires some
regularity on the signal process. 
For instance at the origin we have
\begin{equation}\label{chi-2}
(\tr(P_0)\leq)~r_1\rho(P_0)\leq 1/(4\delta)\Longrightarrow
\EE\left(
\exp{\left[ \delta\Vert X_0-\widehat{X}_0\Vert^{2}\right]}
\right)\leq e.
\end{equation}
The proof of (\ref{chi-2}) and more refined estimates can be found in~\cite{dkt-1-2016}. 

From the numerical viewpoint the trace condition in (\ref{pc}) is related to the initial location of the
particles and the signal-observation perturbations. Signals with a large diffusion part are more 
likely to correct an erroneous initialization. In the same vein, the estimation problems associated with
sensors corrupted by large perturbations are less sensitive to the initialization of the filter.
In the reverse angle, when the signal is almost deterministic and the sensor is highly informative
the particles need to be initialized close to the true value of the signal. 

To better connect our work with existing literature we end our discussion with some connection with 
the variance inflation technique introduced by J.L. Anderson in~\cite{anderson-2,anderson,anderson-3} and further developed
by D.T. B. Kelly, K.J. Law, A. M. Stuart~\cite{kelly} and by X. T. Tong, A. J. Majda and D. Kelly~\cite{tong}. In discrete time settings this technique amounts
of adding an extra positive matrix in the Riccati updating step. This strategy allows to control the fluctuations
of the sample covariance matrices. In continuous time settings, this technique amounts of changing the covariance
matrix $\Pa_{\eta_t}$ in the EKF diffusion (\ref{EEnKF}) by $\Pa_{\eta_t}+\theta~Id$ for some tuning parameter $\theta>0$.
The resulting EKF-diffusion (\ref{EEnKF}) is given by the equation
\begin{eqnarray*}
d\overline{X}_t&=&\left(\Aa\left(\overline{X}_t,\EE[\overline{X}_t~|~\Ga_t]\right)-\theta~S~\overline{X}_t\right)~dt+\Pa_{\eta_t}B^{\prime}R^{-1}_{2}~\left[dY_t-\left(B\overline{X}_t~dt+R^{1/2}_{2}~d\overline{V}_{t}\right)\right]\\
&&+\left[
R^{1/2}_{1}~d\overline{W}_t-\theta~B^{\prime}R^{-1/2}_{2}d\overline{V}_{t}\right]+\theta~B^{\prime}R^{-1}_{2}~dY_t.
\end{eqnarray*}
The stabilizing effects of the variance inflation technique are clear. The last term in the r.h.s. of the above displayed formula has no effect (by simple coupling)
on the stability properties of the diffusion.
The form of the drift also indicates that we increase the Lyapunov exponent by an additional factor $\theta$ (as soon as $\rho(S)>0$). In addition we increase the noise of the diffusion by a factor $\theta^2$, in the sense that the covariance matrix of the perturbation term $R^{1/2}_{1}~d\overline{W}_t-\theta~B^{\prime}R^{-1/2}_{2}d\overline{V}_{t}$ is given
by $R_1+\theta^2S$. We believe that the stability analysis of these regularized models is simplified by these additional regularity properties.
This class of regularized nonlinear diffusions can probably be studied quite easily using the stochastic analysis developed in this article. We plan to develop this analysis in a forthcoming study.
\section{Some preliminary results}\label{preliminary-section}

This short section presents a couple of pivotal results. The first one ensures that
 the Extended Kalman-Bucy filter coincides with the $\Ga_t$-conditional expectations 
of the nonlinear diffusion $\overline{X}_t$. The second result shows that the stochastic processes $\left(m_t,p_t\right)$ satisfy the same equation as
$\left(\widehat{X}_t,P_t\right)$, up to some local fluctuation
orthogonal martingales with angle brackets that only depend on the sample covariance matrix $p_t$.

\begin{prop}\label{prop-EKF-diffusion-mean-cov}
We have the equivalence
$$
\EE(\overline{X}_0)=\widehat{X}_0\quad
\mbox{and}\quad
\Pa_{\eta_0}=P_0
\Longleftrightarrow~\forall t\geq 0\quad
\EE(\overline{X}_t~|~\Ga_t)=\widehat{X}_t\quad\mbox{and}\quad
\Pa_{\eta_t}=P_t .
$$

\end{prop}
\proof
Taking the $\Ga_t$-conditional expectations in (\ref{EEnKF}) we find the diffusion equation
$$
d\,\EE(\overline{X}_t~|~\Ga_t)=A(\EE(\overline{X}_t~|~\Ga_t))~dt+\Pa_{\eta_t}B^{\prime}R^{-1}_{2}~\left[dY_t-B~\EE(\overline{X}_t~|~\Ga_t)dt\right].
$$
Equivalently, if we set $\EE(\overline{X}_t~|~\Ga_t)=\widehat{X}_t$ then we find that
$$
d\,\widehat{X}_t=A(\widehat{X}_t)~dt+\Pa_{\eta_t}B^{\prime}R^{-1}_{2}~\left[dY_t-B~\widehat{X}_t~dt\right].
$$
Let us compute the evolution of $\Pa_{\eta_t}$. We set
$
\widetilde{X}_t=\overline{X}_t-\EE(\overline{X}_t~|~\Ga_t)=\overline{X}_t-\widehat{X}_t
$. In this notation we have
\begin{eqnarray*}
d\widetilde{X}_t
&=&\partial A(\EE(\overline{X}_t~|~\Ga_t))~\widetilde{X}_t~dt~+~R^{1/2}_{1}~d\overline{W}_t-\Pa_{\eta_t}B^{\prime}R^{-1}_{2}~\left[B\widetilde{X}_tdt+R^{1/2}_{2}~d\overline{V}_{t}\right]\\
&=&\left[\partial A(\EE(\overline{X}_t~|~\Ga_t))-\Pa_{\eta_t}S\right]~\widetilde{X}_t~dt~+~R^{1/2}_{1}~d\overline{W}_t-\Pa_{\eta_t}B^{\prime}R^{-1/2}_{2}~d\overline{V}_{t}.
\end{eqnarray*}
This implies that
\begin{eqnarray*}
d(\widetilde{X}_t\widetilde{X}_t^{\prime})
&=&\left\{\left[\partial A(\widehat{X}_t)-P_tS\right]~\widetilde{X}_t\widetilde{X}_t^{\prime}~dt~+
\widetilde{X}_t\widetilde{X}_t^{\prime}\left[\partial A(\widehat{X}_t)-P_tS\right]^{\prime}+(R+\Pa_{\eta_t}S\Pa_{\eta_t})\right\}~dt\\
&& +~\left[R^{1/2}_{1}~d\overline{W}_t-\Pa_{\eta_t}B^{\prime}R^{-1/2}_{2}~d\overline{V}_{t}\right]\widetilde{X}_t^{\prime}+\widetilde{X}_t\left[R^{1/2}_{1}~d\overline{W}_t-\Pa_{\eta_t}B^{\prime}R^{-1/2}_{2}~d\overline{V}_{t}\right]^{\prime}.
\end{eqnarray*}
Taking the $\Ga_t$-conditional expectations we conclude that
\begin{eqnarray*}
\partial_t\Pa_{\eta_t}&=&\left[\partial A(\widehat{X}_t)-\Pa_{\eta_t}S\right]~\Pa_{\eta_t}~dt~+
\Pa_{\eta_t}\left[H(\widehat{X}_t)-\Pa_{\eta_t}S\right]^{\prime}+(R+\Pa_{\eta_t}S\Pa_{\eta_t})\\
&=&\partial A(\widehat{X}_t)\Pa_{\eta_t}+\Pa_{\eta_t}\partial A(\widehat{X}_t)^{\prime}+R-\Pa_{\eta_t}S\Pa_{\eta_t}.
\end{eqnarray*}
This ends the proof of the proposition.
\cqfd

\begin{theo}[Fluctuation theorem~\cite{dt-1-2016}]\label{fluctuation-theorem-ref}
The stochastic processes $\left(m_t,p_t\right)$ defined in (\ref{fv1-3-2}) satisfy the diffusion equations
\begin{equation}\label{f21}
dm_t=A\left[m_t\right]~dt+p_t~B^{\prime}R^{-1}_{2}~\left(dY_t-Bm_t~dt\right)+\frac{1}{\sqrt{N}}~d\overline{M}_t
\end{equation}
with the vector-valued martingale $\overline{M}_t=\left(\overline{M}_t(k)\right)_{1\leq k\leq r_1}$ with the  
angle-brackets
\begin{equation}\label{f22}
\displaystyle
\partial_t\langle \overline{M}_t(k),\overline{M}_t(k^{\prime})\rangle_t=R(k,k^{\prime})+
\left(p_tSp_t\right)(k,k^{\prime}).
\end{equation}
We also have the matrix-valued diffusion
\begin{equation}\label{f23}
dp_t=\left(\partial A\left[m_t\right]p_t+p_t\partial A\left[m_t\right]^{\prime}-p_tSp_t+R\right)~dt+\frac{1}{\sqrt{N-1}}~dM_t
\end{equation}
with a symmetric matrix-valued martingale $M_t=\left(M_t(k,l)\right)_{1\leq k,l\leq r_1}$
and the angle brackets 
\begin{equation}\label{f24}
\begin{array}{rcl}
\displaystyle\partial_t\left\langle  M(k,l), M(k^{\prime},l^{\prime})\right\rangle_t&
=&\left(R+p_tSp_t\right)(k,k^{\prime})~ p_t(l,l^{\prime})+
\left(R+p_tSp_t\right)(l,l^{\prime})~ p_t(k,k^{\prime})\\
&&\\
&&\displaystyle+
\left(R+p_tSp_t\right)(l^{\prime},k) ~p_t(k^{\prime},l)+
\left(R+p_tSp_t\right)(l,k^{\prime}) ~p_t(k,l^{\prime}).
\end{array}
\end{equation}
In addition we have the orthogonality properties
$$ 
\left\langle  M(k,l), \overline{M}(l^{\prime})\right\rangle_t=\left\langle  M(k,l), V(k^{\prime})\right\rangle_t=\left\langle  \overline{M}(l^{\prime}), V(k^{\prime})\right\rangle_t=0
$$
for any $1\leq k,l,l^{\prime}\leq r_1$ and any $1\leq k^{\prime}\leq r_2$.
\end{theo}

\proof
We have
$$
\begin{array}{l}
d(\xi^i_t-m_t)=\left[\partial A\left(m_t\right)-p_tB^{\prime}S\right](\xi^i_t-m_t)dt+dM^i_t
\end{array}
$$
with the martingale
$$
dM^i_t:=R^{1/2}_{1}\left(d\overline{W}_t^i-\frac{1}{N}\sum_{1\leq j\leq N}d\overline{W}_t^j\right)-p_tB^{\prime}R^{-1/2}_{2}
\left(d\overline{V}_t^i-\frac{1}{N}\sum_{1\leq j\leq N}d\overline{V}_t^j\right).
$$
Notice that
$$
\partial_t\langle M^i(k),M^i(k^{\prime})\rangle_t=\left(1-\frac{1}{N}\right)~\left(R+p_tSp_t\right)(k,k^{\prime})
$$
and for $i\not=j$
$$
\partial_t\langle M^i(k),M^j(k^{\prime})\rangle_t=-\frac{1}{N}~\left(R+p_tSp_t\right)(k,k^{\prime}).
$$
The end of the proof follows the proof of~\cite[Theorem 1]{dt-1-2016}, thus it is skipped.
This ends the proof of the theorem.
\cqfd

\section{Stability properties}\label{sec-stability}

This section is dedicated to the long-time behaviour of the EKF-diffusion (\ref{EEnKF}), mainly with the proof of Theorem~\ref{th-intro-1}. We use the stochastic differential inequality calculus developed in~\cite{dt-1-2016,dkt-1-2016}. 
Let $\Ya_t$ be some nonnegative process  defined on 
 some probability space  $(\Omega,\Fa,\PP)$  equipped with a filtration $\Fa=(\Fa_t)_{t\geq 0}$
of $\sigma$-fields. Also let $(\Za_t,\Za^+_t)$  be some processes and $\Ma_t$  be some continuous  $\Fa_t$-martingale. We use the 
following definition
\begin{equation}\label{stoch-ineq}
d\Ya_t\leq \Za_t^+~dt+d\Ma_t~\Longleftrightarrow~ \left(d\Ya_t=\Za_t~dt+d\Ma_t\quad\mbox{\rm with}\quad \Za_t\leq \Za^+_t\right).
\end{equation}
We recall some useful algebraic properties of the above stochastic inequalities.

Let $(\overline{\Ya}_t,\overline{\Za}_t^+,\overline{\Za}_t,\overline{\Ma}_t)$ be another collection of processes satisfying
the above inequalities, and $(\alpha,\overline{\alpha})$ a couple of nonnegative parameters. In this case it is readily checked that
$$
d(\alpha~\Ya_t+\overline{\alpha}~\overline{\Ya}_t)\leq (\alpha~\Za_t^++\overline{\alpha}~\overline{\Za}_t^+)~dt+d(\alpha~\Ma_t+\overline{\alpha}~\overline{\Ma}_t)
$$
and
$$
d(\Ya_t\overline{\Ya}_t)\leq \left[\overline{\Za}_t^+\Ya_t+ \Za_t^+\overline{\Ya}_t+\partial_t\langle \Ma,\overline{\Ma}\rangle_t\right]~dt+\Ya_t~d\overline{\Ma}_t+\overline{\Ya}_t~d\Ma_t.
$$

We consider a couple of diffusions $(\overline{X}_t,\overline{Z}_t)$  coupled with the same Brownian motions $(\overline{V}_t,\overline{W}_t)$ and the same observation process
$Y_t$, and we set
$$
\Fa_t:=\Ga_t\vee\sigma\left((\overline{X}_s,\overline{Z}_s),~s\leq t\right).
$$
Next proposition provides uniform estimates of the $\LL_{\delta}$-centered moments of the EKF-diffusion
with respect to the time horizon. 
\begin{prop}
Assume that $\lambda_{\partial A}>0$. In this situation, for any $\delta\geq 1$ and any time horizon $s\geq 0$ we have the uniform almost sure estimates
\begin{equation}\label{unif-control}
\begin{array}{l}
\EE\left(\Vert \overline{X}_t-\widehat{X}_t\Vert^{\delta}~|~\Fa_s\right)^{2/\delta}
\displaystyle\leq e^{-\lambda_{\partial A}(t-s)}~\Vert\overline{X}_s-\widehat{X}_s\Vert^{2}\\
\\
\hskip3cm+(2\delta-1)\left[\lambda_{R}^{-1}(1+2~\displaystyle(\lambda_{R}\lambda_{S})^{-1})+2
e^{-\lambda_{\partial A}(t+s)}\tr(P_0)^2\lambda_{S}^{-1}
\right].
\end{array}
\end{equation}
\end{prop}

\proof
We have
$$
d(\overline{X}_t-\widehat{X}_t)
=\left[\partial A(\widehat{X}_t)-P_tS\right]~(\overline{X}_t-\widehat{X}_t)~dt~+~R^{1/2}_{1}~d\overline{W}_t-P_tB^{\prime}R^{-1/2}_{2}~d\overline{V}_{t}
$$
This implies that
$$
\begin{array}{l}
d\Vert\overline{X}_t-\widehat{X}_t\Vert^2\\
\\
=\left[2\langle \overline{X}_t-\widehat{X}_t,
\left[\partial A(\widehat{X}_t)-P_tS\right]~(\overline{X}_t-\widehat{X}_t)\rangle+\tr(R_{1})+\tr(P_t^2S)\right]~dt+dM_t\\
\\
\leq \left[-\lambda_{\partial A}~\Vert\overline{X}_t-\widehat{X}_t\Vert^2+\Ua_t\right]~dt+dM_t
\end{array}
$$
with
the process
\begin{eqnarray*}
\Ua_t&:=&\tr(R)+\tr(P_t^2S)\leq \tr(R)+\rho(S)\tr(P_t)^2\\
&\leq&
 \tr(R)+\rho(S)
\left(e^{-\lambda_{\partial A}t}~\tr(P_0)+1/\lambda_{R}\right)^2
\end{eqnarray*}
and the martingale
$$
dM_t:=2\langle \overline{X}_t-\widehat{X}_t,R^{1/2}_{1}~d\overline{W}_t-P_tB^{\prime}R^{-1/2}_{2}~d\overline{V}_{t}
\rangle .
$$
Observe that the angle bracket of this martingale satisfies the property
\begin{eqnarray*}
\partial_t\langle M\rangle_t&=&4\langle  \overline{X}_t-\widehat{X}_t, \left(R+P_tSP_t\right)
(\overline{X}_t-\widehat{X}_t) \rangle\leq 4 \Vert\overline{X}_t-\widehat{X}_t\Vert^2~\Vert R+P_tSP_t\Vert .
\end{eqnarray*}
By~\cite[Corollary 2.2]{dkt-1-2016} for any $\delta\geq 1$ we have
$$
\begin{array}{l}
\EE\left(\Vert \overline{X}_t-\widehat{X}_t\Vert^{\delta}~|~\Fa_s\right)^{2/\delta}
\displaystyle\leq \exp{\left(-\lambda_{\partial A}(t-s)\right)}~\Vert\overline{X}_s-\widehat{X}_s\Vert^{2}\\
\\
\hskip3cm+(2\delta-1)\displaystyle
\int_s^t~\exp{\left(-\lambda_{\partial A}(t-u)\right)}~ \displaystyle
\left(\tr(R)+\rho(S)\tr(P_u)^2\right)~du .
\end{array}$$
Observe that by \eqref{riccati-bound} 
$$
\begin{array}{l}
\displaystyle
\rho(S)\int_s^t~\exp{\left(-\lambda_{\partial A}(t-u)\right)}~ \displaystyle
\tr(P_u)^2~du\\
\\
\leq 2 \rho(S)\displaystyle
\int_s^t~\exp{\left(-\lambda_{\partial A}(t-u)\right)}~ \displaystyle
\left[e^{-2\lambda_{\partial A}u}~\tr(P_0)^2+1/\lambda_{R}^2\right]~du\\
\\
 \leq 2\displaystyle(\lambda_{R}^{2}\lambda_{S})^{-1}+2\exp{\left(-\lambda_{\partial A}(t+s)\right)}\tr(P_0)^2\lambda_{S}^{-1}.
\end{array}$$
This ends the proof of the proposition.
\cqfd

\begin{theo}\label{theo-stability}
When the initial random states $\overline{X}_0$ and $\overline{Z}_0$ have the same
first and second order statistics, that is when 
$
(\widehat{X}_0,P_0)=(\widecheck{X}_0,\widecheck{P}_0)$, we have the almost sure contraction estimates:
$$
\Vert\overline{X}_t-\overline{Z}_t\Vert^2\leq \exp{\left[-\lambda_{\partial A}t\right]}~\Vert\overline{X}_0-\overline{Z}_0\Vert^2.
$$
More generally, when condition (\ref{reference-stability}) is met with $\lambda_S\geq 4^4$,
for any $\epsilon\in [0,1[$ there exists some $s$ such that for any $t\geq s$ and any $1\leq \delta\leq 4^{-4}~\sqrt{\lambda_{S}}$ we have 
\begin{equation}\label{pre-w}
\EE\left(\Vert \overline{X}_t-\overline{Z}_t\Vert^{2\delta}~|~\Fa_{s}\right)^{1/\delta}\leq 
\displaystyle
\exp{\left[-(1-\epsilon)\overline{\lambda}_{\partial A}(t-s)\right]}~\left[
\Vert \overline{X}_{s}-\overline{Z}_{s}\Vert^2+\overline{\Za}_{s}\right]
\end{equation}
with some exponent
$
\overline{\lambda}_{\partial A}\geq \widehat{\lambda}_{\partial A}\wedge({\lambda_{\partial A}}/{2})
$,
and some process  $\overline{\Za}_{t}$ satisfying the uniform moment condition 
\begin{equation}\label{uni-condition-w-check}
\sup_{t\geq 0}{\EE\left( 
\overline{\Za}_{t}^{\alpha/4}
\right)}<\infty
\quad\mbox{ for any}\quad \alpha\leq\lambda_{R,S}~\sqrt{\lambda_{S}}.
\end{equation}

\end{theo}

Before getting into the details of the proof of this theorem we mention that (\ref{w})
is a direct consequence of 
 (\ref{pre-w}) combined with the uniform estimates (\ref{unif-control}). Indeed, applying
 (\ref{pre-w}), for any $\delta\geq 2$ we have
$$
\EE\left(\Vert \overline{X}_t-\overline{Z}_t\Vert^{\delta}\right)^{1/\delta}\leq 
\displaystyle
\exp{\left[-(1-\epsilon)\overline{\lambda}_{\partial A}(t-s)/2\right]}~\left(
\EE\left[
\Vert \overline{X}_{s}-\overline{Z}_{s}\Vert^{\delta}\right]^{1/\delta}
+
\EE\left[
\overline{\Za}_{s}^{\delta/2}
\right]^{1/\delta}\right).
$$
Using  (\ref{unif-control}) and the fact that
$$
1\leq \delta/2\leq  16^{-1}~\sqrt{\lambda_{S}}\leq  {4^{-1}}~\lambda_{R,S}~\sqrt{\lambda_{S}}
$$
we conclude that
$$
\WW_{\delta}(\eta_t,\breve{\eta}_t)\leq c~
\exp{\left[-t~(1-\epsilon)(1-s/t)\overline{\lambda}_{\partial A}/2\right]}\leq
c~
\exp{\left[-t~(1-2\epsilon)\overline{\lambda}_{\partial A}/2\right]}
$$
as soon as $s/t\leq \epsilon$. The end of the proof of (\ref{w}) is now clear.

Now we come to the proof of the theorem.

{\bf Proof of Theorem~\ref{theo-stability}:}

We have
\begin{eqnarray*}
d\overline{X}_t&=&\Aa(\overline{X}_t,\widehat{X}_t)~dt~+~R^{1/2}_{1}~d\overline{W}_t+P_tB^{\prime}R^{-1}_{2}~\left[dY_t-\left(B\overline{X}_tdt+R^{1/2}_{2}~d\overline{V}_{t}\right)\right].
\end{eqnarray*}
Using the decomposition
$$
\begin{array}{rcl}
\widecheck{P}_tS\overline{Z}_t-P_tS\overline{X}_t&=&-P_tS(\overline{X}_t-\overline{Z}_t)+
(\widecheck{P}_t-P_t)S\overline{Z}_t
\end{array}
$$
we readily check that
$$
\begin{array}{l}
d\left(\overline{X}_t-\overline{Z}_t\right)\\
\\
=\left\{\left[\Aa(\overline{X}_t,\widehat{X}_t)-\Aa(\overline{Z}_t,\widecheck{X}_t)\right]-P_tS(\overline{X}_t-\overline{Z}_t)\right\}~dt
+\left[P_t-\widecheck{P}_t\right]S (X_t-\overline{Z}_t)~dt +d\Ma_t
\end{array}
$$
with the martingale
$$
\begin{array}{l}
d\Ma_t:=
\left[P_t-\widecheck{P}_t\right]B^{\prime}R^{-1/2}_{2}~d(V_{t}-\overline{V}_{t})\\
\\
\Rightarrow\begin{array}[t]{rcl}
\partial_t\langle \Ma\rangle_t&=&\Vert \left[P_t-\widecheck{P}_t\right]B^{\prime}R^{-1/2}_{2}\Vert_F^2=\tr\left(\left[P_t-\widecheck{P}_t\right]^2S\right)\leq \Va_t:=\rho(S)~
\Vert P_t-\widecheck{P}_t\Vert_F^2.
\end{array}\end{array}
$$
When the initial random states $\overline{X}_0$ and $\overline{Z}_0$ are possibly different but they have the same
first and second order statistics we have
$$
\widehat{X}_0=\widecheck{X}_0\quad\mbox{\rm and}\quad P_0=\widecheck{P}_0\quad\Longrightarrow~\forall t\geq 0\quad
\widehat{X}_t=\widecheck{X}_t\quad\mbox{\rm and}\quad P_t=\widecheck{P}_t.
$$
In this particular situation we have
$$
\Aa(\overline{X}_t,\widehat{X}_t)-\Aa(\overline{Z}_t,\widecheck{X}_t)=\partial A(
\widecheck{X}_t)~(\overline{X}_t-\overline{Z}_t)
$$
and
$$
\partial_t\left(\overline{X}_t-\overline{Z}_t\right)
=\left[\partial A(
\widecheck{X}_t)-P_tS\right]~(\overline{X}_t-\overline{Z}_t).
$$
This implies that
$$
\partial_t\Vert\overline{X}_t-\overline{Z}_t\Vert^2=2\langle (\overline{X}_t-\overline{Z}_t),
\left[\partial A(
\widecheck{X}_t)-P_tS\right]~(\overline{X}_t-\overline{Z}_t)\rangle\leq -\lambda_{\partial A}~\Vert\overline{X}_t-\overline{Z}_t\Vert^2.
$$
This ends the proof of the first assertion.

More generally, we have
$$
\begin{array}{l}
\Aa(\overline{X}_t,\widehat{X}_t)-\Aa(\overline{Z}_t,\widecheck{X}_t)
=\partial A(
\widecheck{X}_t)~(\overline{X}_t-\overline{Z}_t)\\
\\\hskip2cm+
\left[A(
\widehat{X}_t)-A(
\widecheck{X}_t)\right]-\partial A(
\widecheck{X}_t)(\widehat{X}_t-\widecheck{X}_t)+\left[\partial A(
\widehat{X}_t)-\partial A(
\widecheck{X}_t)\right]
~(\overline{X}_t-\widehat{X}_t).
\end{array}
$$
This yields the  estimate
$$
\begin{array}{l}
\langle \overline{X}_t-\overline{Z}_t,\left(\Aa(\overline{X}_t,\widehat{X}_t)-\Aa(\overline{Z}_t,\widecheck{X}_t)\right)-P_tS (\overline{X}_t-\overline{Z}_t)\rangle\\
\\
\leq \displaystyle-\frac{\lambda_{\partial A}}{2}~\Vert\overline{X}_t-\overline{Z}_t\Vert^2+\langle \overline{X}_t-\overline{Z}_t,\left[\partial A(
\widehat{X}_t)-\partial A(
\widecheck{X}_t)\right]
~(\overline{X}_t-\widehat{X}_t)
\rangle\\
\\
\hskip3cm+\langle \overline{X}_t-\overline{Z}_t,\left[A(
\widehat{X}_t)-A(
\widecheck{X}_t)\right]-\partial A(
\widecheck{X}_t)(\widehat{X}_t-\widecheck{X}_t)\rangle\\
\\
\leq \displaystyle-\frac{\lambda_{\partial A}}{2}~\Vert\overline{X}_t-\overline{Z}_t\Vert^2+~\Vert
\widehat{X}_t-\widecheck{X}_t\Vert~\Vert\overline{X}_t-\overline{Z}_t\Vert\left(\kappa_{\partial A}
~\Vert\overline{X}_t-\widehat{X}_t\Vert+\kappa_{\partial A}+\Vert\partial A\Vert\right)
\end{array}
$$
We also have
$$
\begin{array}{l}
\langle \overline{X}_t-\overline{Z}_t,\left[P_t-\widecheck{P}_t\right]S (X_t-\overline{Z}_t)\rangle\leq 
\Vert P_t-\widecheck{P}_t\Vert_F~\Vert \overline{X}_t-\overline{Z}_t\Vert~\Vert S (X_t-\overline{Z}_t)\Vert .
\end{array}
$$
This implies that
$$
\begin{array}{l}
d\Vert\overline{X}_t-\overline{Z}_t\Vert^2\\
\\
\leq\displaystyle\left[-\lambda_{\partial A}~\Vert\overline{X}_t-\overline{Z}_t\Vert^2+2~\Vert
\widehat{X}_t-\widecheck{X}_t\Vert~\Vert\overline{X}_t-\overline{Z}_t\Vert\left(\kappa_{\partial A}
~\Vert\overline{X}_t-\widehat{X}_t\Vert+\kappa_{\partial A}+\Vert\partial A\Vert\right)\right]~dt\\
\\
+\left[2\Vert P_t-\widecheck{P}_t\Vert_F~\Vert \overline{X}_t-\overline{Z}_t\Vert~\Vert S (X_t-\overline{Z}_t)\Vert\right]~dt
 +2\sqrt{\Va_t}~\Vert \overline{X}_t-\overline{Z}_t\Vert~d\overline{\Ma}_t
\end{array}
$$
with $$
\Va_t=\rho(S)~
\Vert P_t-\widecheck{P}_t\Vert_F^2$$
and a rescaled continuous martingale $\overline{\Ma}_t$ such that
$
\partial_t\langle \overline{\Ma}\rangle_t\leq  1$. On the other hand, we have
$$
\begin{array}{l}
2~\Vert\overline{X}_t-\overline{Z}_t\Vert~\Vert
\widehat{X}_t-\widecheck{X}_t\Vert\left(\kappa_{\partial A}
~\Vert\overline{X}_t-\widehat{X}_t\Vert+\kappa_{\partial A}+\Vert\partial A\Vert\right)\\
\\
\leq \displaystyle\frac{\lambda_{\partial A}}{4}\Vert\overline{X}_t-\overline{Z}_t\Vert^2+\frac{4}{\lambda_{\partial A}}\Vert\widehat{X}_t-\widecheck{X}_t\Vert^2\left(\kappa_{\partial A}
~\Vert\overline{X}_t-\widehat{X}_t\Vert+\kappa_{\partial A}+\Vert\partial A\Vert\right)^2
\end{array}
$$
and
$$
\begin{array}{l}
2\Vert P_t-\widecheck{P}_t\Vert_F~\Vert \overline{X}_t-\overline{Z}_t\Vert~\Vert S (X_t-\overline{Z}_t)\Vert\\
\\
\leq\displaystyle \frac{\lambda_{\partial A}}{4}\Vert \overline{X}_t-\overline{Z}_t\Vert^2+\frac{4}{\lambda_{\partial A}}~
\Vert P_t-\widecheck{P}_t\Vert_F^2~\Vert S (X_t-\overline{Z}_t)\Vert^2.
\end{array}
$$
We conclude that
$$
d\Vert \overline{X}_t-\overline{Z}_t\Vert^2\leq  \left[-\frac{\lambda_{\partial A}}{2}~\Vert\overline{X}_t-\overline{Z}_t\Vert^2+\Ua_t\right]~dt+2\sqrt{\Va_t}~\Vert \overline{X}_t-\overline{Z}_t\Vert~d\overline{\Ma}_t
$$
with 
$$
\Ua_t:=\alpha_t~\Vert
\widehat{X}_t-\widecheck{X}_t\Vert^2+\beta_t~\Vert P_t-\widecheck{P}_t\Vert_F^2
$$
and the parameters
$$
\alpha_t:=\frac{4}{\lambda_{\partial A}}\left(\kappa_{\partial A}
~\Vert\overline{X}_t-\widehat{X}_t\Vert+\kappa_{\partial A}+\Vert\partial A\Vert\right)^2\quad\mbox{\rm and}\quad\beta_t:=\frac{4}{\lambda_{\partial A}}~\Vert S (X_t-\overline{Z}_t)\Vert^2.
$$
By (\ref{reference-stability}) and (\ref{unif-control}), for any $\delta\leq 2^{-1}\sqrt{\lambda_{S}}$ and any $t\geq s$ we  have
\begin{eqnarray*}
\EE\left(\alpha_t^{\delta/4}~\Vert
\widehat{X}_t-\widecheck{X}_t\Vert^{\delta/2}~|~\Fa_{s}\right)^{4/\delta}&\leq&
\EE\left(\Vert
\widehat{X}_t-\widecheck{X}_t\Vert^{\delta}~|~\Fa_{s}\right)^{2/\delta}~\EE\left(\alpha_t^{\delta/2}~|~\Fa_{s}\right)^{2/\delta}\\
&\leq &  \overline{\Za}_{s}~\displaystyle 
\exp{\left(-\widehat{\lambda}_{\partial A} (1-\epsilon) (t-s)\right)}~
\end{eqnarray*}
for  some process $\overline{\Za}_{s}$ satisfying the uniform moment condition (\ref{uni-condition-w-check}).
In the same vein we check that
$$
\EE\left(\Ua_t^{\delta/4}~|~\Fa_{s}\right)^{4/\delta}\vee \EE\left(\Va_t^{\delta/4}~|~\Fa_{s}\right)^{4/\delta}\leq
\overline{\Za}_{s}~\displaystyle 
\exp{\left(-\widehat{\lambda}_{\partial A} (1-\epsilon) (t-s)\right)}
$$
for any $s\geq t_0$.
By~\cite[Corollary 2.2]{dkt-1-2016} we have
$$
\begin{array}{l}
\EE\left(\Vert \overline{X}_t-\overline{Z}_t\Vert^{\delta/4}~|~\Fa_{s}\right)^{8/\delta}\\
\\
\displaystyle\leq \displaystyle \exp{\left(-\left[
\frac{\lambda_{\partial A}}{2}(t-s)\right]\right)}~\Vert \overline{X}_{s}-\overline{Z}_{s}\Vert^{2}\\
\\
\hskip4cm\displaystyle+~n~\overline{\Za}_{s}
\int_{s}^t~\exp{\left(-\left[
\frac{\lambda_{\partial A}}{2}(t-u)+\widehat{\lambda}_{\partial A} (1-\epsilon) (u-s)\right]\right)}~du\\
\\
 \displaystyle\leq \displaystyle e^{-
\frac{\lambda_{\partial A}}{2}(t-s)}~\Vert \overline{X}_{s}-\overline{Z}_{s}\Vert^{2}+\displaystyle~\frac{n~\overline{\Za}_{t_0}}{\vert \widehat{\lambda}_{\partial A} (1-\epsilon)-{\lambda_{\partial A}}/{2}\vert}~
\vert e^{-\frac{\lambda_{\partial A}}{2}(t-s)}-e^{-\widehat{\lambda}_{\partial A} (1-\epsilon) (t-s)}\vert.
\end{array}
$$
The end of the proof of the theorem is now easily completed.
\cqfd

\section{Quantitative propagation of chaos estimates}\label{pc-section}

\subsection{Laplace exponential moment estimates}\label{Laplace-section}
The analysis of EKF filters and their particle interpretation is mainly based on the estimation of the stochastic
exponential function 
$$
\Ea_{\Gamma}(t):=\exp{\left[\int_0^t\Gamma_A(s)~ ds\right]}
$$
with the stochastic functional
$$
\Gamma_A(s):=-
\left[
\lambda_{\partial A}-\left(2\kappa_{\partial A}~\tr(P_t)+\rho(S)~\Vert X_t- \widehat{X}_t\Vert~\right)\right].
$$
Assume condition  (\ref{reference-stability})  is satisfied and set
$$
\Lambda_{\partial A}\left[\epsilon,\delta\right]/\lambda_{\partial A}:=
1-\frac{2}{\lambda_{K}\lambda_R}+\frac{1}{\lambda_{S}}\left(\frac{3}{4}-\delta\right)-
\frac{1}{\delta}~
\frac{\epsilon\lambda_A}{2\lambda_{\partial A}}.
$$
Observe that for any $\delta> 0$ we have
$$
\epsilon=\frac{1}{2}~\frac{\lambda_{\partial A}}{\lambda_A}\Longrightarrow
\Lambda_{\partial A}\left[\epsilon,\sqrt{\lambda_S}/2\right]=
\widehat{\lambda}_{\partial A}\geq \Lambda_{\partial A}\left[\epsilon,\delta\right].
$$

The next technical lemma provides some key $\delta$-exponential moments estimates. Its proof is quite technical, thus it is housed in the appendix, Section~\ref{Proof-Laplace-Sec}.

\begin{lem}\label{lem-Laplace}
\begin{itemize}
\item For any $\delta> 0$ and any $0\leq s\leq t$ we have the almost sure estimate 
\begin{equation}\label{estimation-neg}
\EE\left((\Ea_{\Gamma}(t)/\Ea_{\Gamma}(s))^{-\delta}~|~\Fa_s\right)^{1/\delta}\leq \exp{\left(~\Lambda^-_{\Gamma}~(t-s)\right)}
\quad\mbox{with}
\quad \Lambda^-_{\Gamma}=
\lambda_{\partial A}\left[1-\frac{2}{\lambda_K\lambda_{R}}\right].
\end{equation}
\item For any $\epsilon \in [0,1]$, any  $0<\delta\leq e~\epsilon~\lambda_{R,S}$ and any initial covariance matrix $P_0$ such that
$$
 \tr(P_0)^2\leq \sigma^2(\epsilon,\delta):=\frac{\lambda_S}{\lambda_R}
 \left[\frac{1}{2}+\frac{1}{\lambda_{R}\lambda_S}\right]\left(
e~\epsilon{\lambda_{R,S}}/{\delta}-1\right)
$$
for any time horizon $t\geq 0$ we have the exponential $\delta$-moment estimate
\begin{equation}\label{estimation-positive}
\EE\left[\Ea_{\Gamma}(t)^{\delta}\right]^{1/\delta}\leq c_{\delta}(P_0)~ \exp{\left[\Lambda^+_{\Gamma}(\epsilon,\delta)~t\right]}
\end{equation}
with the parameters
\begin{eqnarray*}
\Lambda^+_{\Gamma}(\epsilon,\delta)&:=& 2\kappa_{\partial A}\sigma(\epsilon,\delta)-\Lambda_{\partial A}\left[\epsilon,\delta\right]-\left(\delta-1\right)\rho(S)
\\
c_{\delta}(P_0)&:=&\exp{\left(1/\delta+\delta\chi(P_0)/(2\lambda_{S})^2\right)}.
\end{eqnarray*}
\item For any $\epsilon\in ]0,1]$  there exists some time horizon $s$ such that
for any $t\geq s$ and any $\delta\leq \sqrt{\lambda_S}/2$ we have the almost sure estimate
\begin{equation}\label{estimation-lemma}
\displaystyle\EE\left(
\Ea_{\Gamma}(t)^{\delta}~|~\Fa_s\right)^{1/\delta}\leq \displaystyle \Ea_{\Gamma}(s)~ \Za_{s}~\displaystyle \exp{\left(-\left\{(1-\epsilon)\widehat{\lambda}_{\partial A}+(\delta-1)\rho(S)\right\}~(t-s)\right)
}
\end{equation}
for some positive random process $\Za_{t}$ s.t. 
$$
\forall \alpha\leq \lambda_{R,S}~\sqrt{\lambda_S},~\qquad
\sup_{t\geq 0}{\EE\left(\Za_t^{\alpha}\right)}<\infty .
$$
\end{itemize}
\end{lem}

\subsection{A non asymptotic convergence theorem}

This section is mainly concerned with the estimation of the $\delta$-moments
of the square errors
$$
\Xi_t:=\Vert(m_t,p_t)-(\widehat{X}_t,P_t)\Vert^2=\Vert m_t-\widehat{X}_t\Vert^2+\Vert p_t-P_t\Vert_F^2.
$$
The analysis is based on a couple of technical lemmas. 

The first one provides uniform moments
estimates with respect to the time parameter. 
\begin{lem}\label{prop-unif-control}
There exists some $\nu> 0$ such that
for any $1\leq n\leq 1+\nu N$ we have
$$
\sup_{t\geq 0}{\EE\left(\tr(p_t)^n\right)}<\infty ,
\qquad
\sup_{t\geq 0}{\EE\left(\Vert \xi^1_t\Vert^n\right)}<\infty\quad
\mbox{and}\quad
\sup_{t\geq 0}{\EE\left(\Vert \zeta^1_t\Vert^n\right)}<\infty .
$$
\end{lem}

The second technical lemma provides a differential perturbation inequality in terms
of the Laplace functionals discussed in Section~\ref{Laplace-section}.

\begin{lem}\label{lem-perturbation}
We have the stochastic differential inequality
\begin{eqnarray*}
d\Xi_t
&\leq &\Xi_t~\left[\Gamma_A(t)+
\sqrt{2\rho(S)}~d\,\Upsilon_t^{(1)}\right]+\left[\Va_t~dt+\sqrt{\Va_t~\Xi_t}~d\,\Upsilon_t^{(2)}\right]
\end{eqnarray*}
with a couple of orthogonal martingales s.t. $\partial_t\langle \Upsilon_{\cdot}^{(i)},\Upsilon_{\cdot}^{(j)}\rangle_t\leq 1_{i=j}$ and some nonnegative 
process $\Va_t$ such that
$$
\sup_{t\geq 0}{\EE\left(\Va^n_t\right)}^{1/n}\leq c(n)/N
\quad\mbox{
for any $1\leq n\leq 1+\nu N$ and some $\nu>0$.}
$$
\end{lem}

The proofs of these two lemmas are rather technical thus they are provided in the appendix, Section~\ref{Lemma-Unif-control-sec} and 
Section~\ref{Proof-perturbation-lemma}.
We are now in position to state and to prove the main result of this section.
\begin{theo}\label{theo-c}
Assume that $(2^{-1}\sqrt{\lambda_S})\wedge (e\lambda_{R,S})\geq 2$. In this situation,
there exist some $N_0\geq 1$ and some $\alpha\in ]0,1]$ such that for any $N_0\leq N$, 
$1\leq \delta\leq (4^{-1}\sqrt{\lambda_S})\wedge (2^{-1}e\lambda_{R,S})$ and any
initial covariance matrix $P_0$ of the signal  we have
 the uniform estimates
$$
 \tr(P_0)^2\leq\frac{1}{2}~\frac{\lambda_S}{\lambda_R}
 \left[1+\frac{2}{\lambda_{R}\lambda_S}\right]\Longrightarrow
\sup_{t\geq 0}{\EE[\Xi_{t}^{\delta}]^{1/\delta}}\leq c/N^{\alpha}.
$$
\end{theo}
\proof
We set
$$
\overline{\Ea}(t):=\Ea_{\Gamma}(t)\Ea_{\Upsilon}(t)=e^{\La_t}
$$
with the exponential martingale
$$
\Ea_{\Upsilon}(t):=\exp{\left[\sqrt{2\rho(S)}~\Upsilon_t^{(1)}-\rho(S) t\right]}$$
and the stochastic process
$$
\La_t:=\int_0^t\Gamma_A(u)~ du+\sqrt{2\rho(S)}~\Upsilon_t^{(1)}-\rho(S) t .
$$

Observe that for any $\delta\geq 0$ we have
\begin{eqnarray*}
\Ea_{\Upsilon}^{-\delta}(t)&=&\exp{\left[-\delta \sqrt{2\rho(S)}~\Upsilon_t^{(1)}+\delta\rho(S) t\right]}=\exp{\left[\delta(1+2\delta)\rho(S) t\right]}~~\Ea_{-2\delta\Upsilon}^{1/2}(t)
\end{eqnarray*}
with the exponential martingale
$$
\Ea_{-2\delta\Upsilon}(t):=\exp{\left[-2\delta \sqrt{2\rho(S)}~\Upsilon_t^{(1)}-4\delta^2\rho(S) t\right]}.
$$
In the same vein we have
\begin{eqnarray*}
\Ea_{\Upsilon}^{\delta}(t)&=&\exp{\left[\delta \sqrt{2\rho(S)}~\Upsilon_t^{(1)}-\delta\rho(S) t\right]}=\exp{\left[\delta(2\delta-1)\rho(S) t\right]}~~\Ea_{2\delta\Upsilon}^{1/2}(t)
\end{eqnarray*}
with the exponential martingale
$$
\Ea_{2\delta\Upsilon}(t):=\exp{\left[2\delta \sqrt{2\rho(S)}~\Upsilon_t^{(1)}-4\delta^2\rho(S) t\right]}.
$$

This yields the estimates
\begin{eqnarray*}
\EE\left(\overline{\Ea}^{-\delta}(t)\right)&=&\exp{\left(\delta(1+2\delta)\rho(S) t\right)}~\EE\left[\Ea_{\Gamma}(t)^{-\delta}~\Ea_{-2\delta\Upsilon}^{1/2}(t)\right]\\
&\leq &\EE\left[\Ea_{\Gamma}(t)^{-2\delta}\right]^{1/2}~
\exp{\left(\delta(1+2\delta)~\rho(S) t\right)}
\\
\EE\left(\overline{\Ea}^{\delta}(t)\right)
&\leq &\EE\left[\Ea_{\Gamma}(t)^{2\delta}\right]^{1/2}~
\exp{\left(\delta(2\delta-1)~\rho(S) t\right)}.
\end{eqnarray*}

Using (\ref{estimation-neg}) and (\ref{estimation-positive}) we find the estimates
\begin{eqnarray}
\EE\left(\overline{\Ea}^{-\delta}(t)\right)^{1/\delta}
&\leq &
\exp{\left(\left[(1+2\delta)~\rho(S)+\Lambda_{\Gamma}^-\right] t\right)} \label{estimation-neg-proof}\\
\EE\left(\overline{\Ea}^{\delta}(t)\right)^{1/\delta}
&\leq &c_{\delta}(P_0)~
\exp{\left(\left[(2\delta-1)~\rho(S)+\Lambda^+_{\Gamma}(\epsilon,\delta)\right] t\right)}. \label{estimation-pos-proof}
\end{eqnarray}
The estimate (\ref{estimation-pos-proof}) is valid for any $\epsilon \in [0,1]$  and any
$$
\displaystyle\delta\leq e~\epsilon~\lambda_{R,S}\quad\mbox{and}\quad
 \tr(P_0)\leq \sigma(\epsilon,\delta).
$$

Using the fact that
\begin{eqnarray*}
d\overline{\Ea}^{-1}(t)&\leq &-e^{-\La_t}~\left(\Gamma_A(t)~dt+\sqrt{2\rho(S)}~d\,\Upsilon_t^{(1)}-\rho(S) dt\right)+\frac{1}{2}~e^{-\La_t}~2\rho(S)~\partial_t\langle \Upsilon^{(1)}\rangle_t~dt\\
&\leq &-\overline{\Ea}^{-1}(t)~\left(\Gamma_A(t)~dt+\sqrt{2\rho(S)}~d\,\Upsilon_t^{(1)}\right)
\end{eqnarray*}
we find the stochastic inequality
\begin{eqnarray*}
d(\Xi_t~\overline{\Ea}^{-1}(t))&\leq &\overline{\Ea}^{-1}(t)~d\Xi_t+\Xi_t ~d\overline{\Ea}^{-1}(t)-2\overline{\Ea}^{-1}(t)~\Xi_t~\rho(S)~dt\\
&\leq &\overline{\Ea}^{-1}(t)~\Xi_t~\left[\Gamma_A(t)+
\sqrt{2\rho(S)}~d\,\Upsilon_t^{(1)}\right]+\overline{\Ea}^{-1}(t)\left[\Va_t~dt+\sqrt{\Va_t~\Xi_t}~d\,\Upsilon_t^{(2)}\right]\\
&&\hskip.3cm-\overline{\Ea}^{-1}(t)~\Xi_t~\left[\Gamma_A(t)~dt+\sqrt{2\rho(S)}~d\,\Upsilon_t^{(1)}\right]-2\overline{\Ea}^{-1}(t)~\Xi_t~\rho(S)~dt\\
&=&\overline{\Ea}^{-1}(t)\left[\left(\Va_t-2~\Xi_t~\rho(S)\right)~dt+\sqrt{\Va_t~\Xi_t}~d\,\Upsilon_t^{(2)}\right].
\end{eqnarray*}
For any $\delta\geq 2$, this implies that
\begin{eqnarray*}
d(\Xi_t~\overline{\Ea}^{-1}(t))^{\delta}&\leq&\delta~\Xi_t^{\delta-1}~\overline{\Ea}^{-\delta}(t)~\left[\left(\Va_t-2~\Xi_t~\rho(S)\right)~dt+\sqrt{\Va_t~\Xi_t}~d\,\Upsilon_t^{(2)}\right]\\
&&+
\delta~\Xi_t^{\delta-1}~\overline{\Ea}(t)^{-\delta}~\frac{(\delta-1)}{2}\Va_t~dt\\
&=&\delta~\Xi_t^{\delta-1}~\overline{\Ea}^{-\delta}(t)~\left[\left(\frac{\delta+1}{2}~\Va_t-2~\Xi_t~\rho(S)\right)~dt+\sqrt{\Va_t~\Xi_t}~d\,\Upsilon_t^{(2)}\right].
\end{eqnarray*}
Taking the expectation we obtain
\begin{eqnarray*}
\partial_t\EE\left[(\Xi_t~\overline{\Ea}^{-1}(t))^{\delta}\right]&\leq&
\frac{\delta(\delta+1)}{2}~
\EE\left[~\left(\Xi_t~\overline{\Ea}^{-1}(t)\right)^{\delta-1}~\overline{\Ea}^{-1}(t) ~\Va_t\right]\\
&&
~-2\delta~\rho(S)~\EE\left[\left(\Xi_t~\overline{\Ea}^{-1}(t)\right)^{\delta}\right].
\end{eqnarray*}
On the other hand using Lemma~\ref{lem-perturbation} and the Laplace estimate (\ref{estimation-neg-proof}) we have
%$$
%\begin{array}{l}
\begin{eqnarray*}
\EE\left(\left(\Xi_t~\overline{\Ea}^{-1}(t)\right)^{\delta-1}~\overline{\Ea}^{-1}(t) ~\Va_t\right)
&\leq &
\EE\left(\left(\Xi_t~\overline{\Ea}^{-1}(t)\right)^{\delta}\right)^{1-1/\delta}
\EE\left(\overline{\Ea}^{-2\delta}(t)\right)^{1/(2\delta)}
\EE\left(\Va_t^{2\delta}\right)^{1/(2\delta)}\\
%\\
%\displaystyle
&\leq & \frac{c}{N}~\exp{\left(\left[(1+4\delta)~\rho(S)+\Lambda_{\Gamma}^-\right] t\right)}~
\EE\left(\left(\Xi_t~\overline{\Ea}^{-1}(t)\right)^{\delta}\right)^{1-1/\delta}.
%\end{array}
%$$
\end{eqnarray*}
This yields
\begin{eqnarray*}
\partial_t\EE\left((\Xi_t\overline{\Ea}^{-1}(t))^{\delta}\right)^{1/\delta}
&\leq & \frac{1}{\delta}~\EE\left((\Xi_t\overline{\Ea}^{-1}(t))^{\delta}\right)^{\frac{1}{\delta}-1}~\partial_t\EE\left((\Xi_t\overline{\Ea}^{-1}(t))^{\delta}\right)\\
&\leq & -2\rho(S)\EE\left((\Xi_t\overline{\Ea}^{-1}(t))^{\delta}\right)^{1/\delta}+\frac{(\delta+1)}{2} \frac{c}{N}~\exp{\left(\left[(1+4\delta)\rho(S)+\Lambda_{\Gamma}^-\right] t\right)} 
\end{eqnarray*}
from which we conclude that
\begin{eqnarray*}
\EE\left((\Xi_t~\overline{\Ea}^{-1}(t))^{\delta}\right)^{1/\delta}&\leq& \exp{\left\{-2\rho(S)t\right\}}~\EE\left(\Xi_0^{\delta}\right)^{1/\delta}+\frac{c}{N}~\exp{\left\{\left((1+4\delta)\rho(S)+\Lambda_{\Gamma}^-\right) t\right\}} \\
&\leq & \frac{c}{N}~\exp{\left\{\left((1+4\delta)\rho(S)+\Lambda_{\Gamma}^-\right) t\right\}} .
\end{eqnarray*}
By Cauchy-Schwarz inequality we also have
$$
\EE\left(\Xi_t^{\delta/2}\right)^{2/\delta}=\EE\left(\overline{\Ea}(t)^{\delta/2}~\left(\Xi_t~\overline{\Ea}^{-1}(t)\right)^{\delta/2}\right)^{2/\delta}\leq \EE\left((\Xi_t~\overline{\Ea}^{-1}(t))^{\delta}\right)^{1/\delta}
\EE\left(\overline{\Ea}^{\delta}(t)\right)^{1/\delta}.
$$
Using (\ref{estimation-pos-proof}) we conclude that
for any $\epsilon \in [0,1]$  and any
$
\displaystyle\delta\leq e~\epsilon~\lambda_{R,S}$ and $
 \tr(P_0)\leq \sigma(\epsilon,\delta),
$

\begin{equation}\label{first-estimate}
\EE\left(\Xi_t^{\delta/2}\right)^{2/\delta}
\leq c_{\delta}(P_0)~\frac{c}{N}~\exp{\left\{\left(6\delta\rho(S)+\Lambda_{\Gamma}^-+\Lambda^+_{\Gamma}(\epsilon,\delta)\right) t\right\}} . 
\end{equation}

On the other hand, by~\cite[Theorem 2.1]{dkt-1-2016} for any $\delta\geq 1$ we also have
\begin{equation}\label{gronwall-n+}
\begin{array}{l}
\EE\left(\Xi_t^{\delta/2}~|~\Fa_s\right)^{2/\delta}
\displaystyle\leq \EE\left[\exp{\left(\delta~\int_s^t\left\{\Gamma_A(u)~+(\delta-1)\rho(S)\right\}~du\right)}~|~\Fa_s\right]^{1/\delta}\\
\\
\hskip3cm\times\displaystyle\left\{
\Xi_s+\frac{1}{N}~\frac{\delta+1}{2}~\int_s^t~\EE\left[\overline{\Va}_u^{\delta}~|~\Fa_s\right]^{1/\delta}~du
\right\}
\end{array}\end{equation}
with the rescaled process 
$$
\overline{\Va}_t:=\exp{\left(\int_s^t\left[-\Gamma_A(u)+2(1-\delta)\rho(S)\right]du\right)}~\Va_t
$$
of the process $\Va_t$ defined in Lemma~\ref{lem-perturbation}.

On the other hand using (\ref{estimation-lemma}) for any $\epsilon\in ]0,1]$  there exists some time horizon $s=s(\epsilon)$ such that
for any $t\geq s$ and any $\delta\leq \frac{1}{2}~\sqrt{\lambda_S}$ we have the almost sure estimate
$$
\EE\left(\Xi_t^{\delta/2}~|~\Fa_s\right)^{2/\delta}
\displaystyle\leq \displaystyle \Za_{s}~\displaystyle \exp{\left(-\left(1-\epsilon\right)\widehat{\lambda}_{\partial A}(t-s)\right)}\times\displaystyle\left\{
\Xi_s+\frac{\delta+1}{2}~\int_s^t~\EE\left[\overline{\Va}_u^{\delta}~|~\Fa_s\right]^{1/\delta}~du
\right\}
$$
with some process $\Za_{s}$ such that
$$
\sup_{t\geq 0}{\EE\left( 
\Za^{\alpha}_t
\right)}<\infty
\quad\mbox{ for any}\quad \alpha\leq \frac{2}{e}~\sqrt{\lambda_{S}} \left(\leq \frac{1}{2}~\lambda_{R,S}~\sqrt{\lambda_{S}}\right).
$$
Combining Cauchy-Schwarz inequality with (\ref{estimation-neg}) and
Lemma~\ref{lem-perturbation} we readily check that
$$
\begin{array}{l}
\displaystyle\EE\left[\overline{\Va}_u^{\delta}~|~\Fa_s\right]^{1/\delta}\\
\\
= \displaystyle\EE\left[{\Va}_u^{\delta}~\exp{\left(\delta\int_s^u\left[-\Gamma_A(v)+2(1-\delta)\rho(S)\right]dv\right)}~|~\Fa_s\right]^{1/\delta}\\
 \\
 \displaystyle\leq 
  \EE\left[\Va_u^{2\delta}\right]^{1/(2\delta)}~\exp{\left(2(1-\delta)\rho(S)(u-s)\right)} ~\displaystyle\EE\left[(\Ea_{\Gamma}(u)/\Ea_{\Gamma}(s))^{-2\delta}~|~\Fa_s\right]^{1/(2\delta)}
\\
\\
\displaystyle\leq 
\frac{c}{N}~\exp{\left[\left(2(1-\delta)\rho(S)+\Lambda^-_{\Gamma}\right)(u-s)\right]} .
\end{array}$$

This yields the estimate
$$
\begin{array}{l}
\EE\left(\Xi_t^{\delta/2}~|~\Fa_s\right)
\displaystyle\leq \displaystyle \Za_{s}^{\delta/2}~\displaystyle \exp{\left(-\frac{\delta}{2}~\left(1-\epsilon\right)~\widehat{\lambda}_{\partial A}~(t-s)\right)}\\
\\
\hskip3cm\times\displaystyle\left\{
\Xi_s+\frac{c}{N}~\exp{\left[\left(2(1-\delta)\rho(S)+\Lambda^-_{\Gamma}\right)(t-s)\right]}
\right\}^{\delta/2} .
\end{array}
$$
This implies that for any $1\leq \delta/2\leq \frac{1}{4}~\sqrt{\lambda_S}$ we have
$$
\begin{array}{l}
\EE\left(\Xi_{t}^{\delta/2}~|~\Fa_s\right)
\displaystyle\leq \displaystyle ~\displaystyle c~\Za_{s}^{\delta/2}~\exp{\left(-\frac{\delta}{2}~\left(1-\epsilon\right)~\widehat{\lambda}_{\partial A}~(t-s)\right)}\\
\\
\hskip3cm\times\displaystyle\left\{
\Xi_s^{\delta/2}+\frac{1}{N^{\delta/2}}~\exp{\left[\frac{\delta}{2}~\left(2(1-\delta)\rho(S)+\Lambda^-_{\Gamma}\right)(t-s)\right]}
\right\}.
\end{array}
$$
Taking the expectation and choosing $\epsilon\leq 1/2$, there exists some time horizon $t_0$ such that
for any $s\geq 0$ and any $\tau\geq s+t_0$
$$
\EE\left(\Xi_{\tau}^{\delta/2}\right)^{2/\delta} \leq  c~\exp{\left(-\widehat{\lambda}_{\partial A}~(\tau-(s+t_0))/2\right)}\displaystyle\left\{1+\frac{1}{N}~\exp{\left[\left(2(1-\delta)\rho(S)+\Lambda^-_{\Gamma}\right)(\tau-(s+t_0)\right]}
\right\}
$$
for any $2\leq \delta\leq 1+\nu N$ for some $\nu>0$, and for some finite constant $c(\delta)<\infty$.
This implies that for any time horizon $t\geq 0$ and any $$2\leq \delta\leq 2^{-1}\sqrt{\lambda_S}\wedge (1+\nu N)$$ we have 
$$
\EE\left(\Xi_{s+t_0+t}^{\delta/2}\right)^{2/\delta}
\displaystyle\leq \displaystyle \displaystyle c\exp{\left(-\frac{\widehat{\lambda}_{\partial A}}{2}~t\right)}\displaystyle\left\{1+\frac{1}{N}~\exp{\left[\left(2(1-\delta)\rho(S)+\Lambda^-_{\Gamma}\right)t\right]}
\right\}.
$$
This yields the uniform estimates
$$
\displaystyle\sup_{u\in [t+t_0,\infty[}{\EE\left(\Xi_{u}^{\delta/2}\right)^{2/\delta}}
=\sup_{s\geq 0}{\EE\left(\Xi_{s+t_0+t}^{\delta/2}\right)^{2/\delta}}
\displaystyle\leq  \displaystyle \displaystyle c~\left\{\exp{\left[-\frac{\widehat{\lambda}_{\partial A}}{2}~t\right]}\displaystyle+
\frac{1}{N}~\exp{\left[\lambda_{\Gamma}t\right]}\right\}
$$
with the parameters
\begin{eqnarray*}
\lambda_{\Gamma}&:=&
\Lambda^-_{\Gamma}-2\rho(S)=
\frac{\lambda_{\partial A}}{2}\left[\left(1-\frac{4}{\lambda_K\lambda_{R}}\right)+\left(1-\frac{4}{\lambda_S}\right)\right]>0.
\end{eqnarray*}

On the other hand, by (\ref{first-estimate}) for any time horizon $t\geq 0$ and any
$
\displaystyle\delta\leq e~\lambda_{R,S}$ and any $P_0$ s.t. $
 \tr(P_0)\leq \sigma(1,e\lambda_{R,S}/2)
$ we have the uniform estimates
\begin{eqnarray*}
\sup_{s\in [0,t_0+t]}\EE\left(\Xi_s^{\delta/2}\right)^{2/\delta}
&\leq &c~\frac{1}{N}~\exp{\left[\lambda_{\Gamma}^{\prime} t\right]}  
\end{eqnarray*}
with
\begin{eqnarray*}
\lambda_{\Gamma}^{\prime}&:=& 5e\lambda_{\partial A}~\lambda_{R,S}/\lambda_S+\Lambda_{\Gamma}^-+\Lambda^+_{\Gamma}(1,e\lambda_{R,S}/2).
\end{eqnarray*}
We conclude that for any time horizon $t\geq 0$ 
$$
\displaystyle\sup_{s\geq 0}{\EE\left(\Xi_{s}^{\delta/2}\right)^{2/\delta}}\leq c~\left\{
~\exp{\left[-\frac{\widehat{\lambda}_{\partial A}}{2}~t\right]}\displaystyle+ 
\frac{1}{N}~\exp{\left[(\lambda_{\Gamma}\vee \lambda_{\Gamma}^{\prime})t\right]}\right\}.
$$
Choosing $t=t(N)$ such that
$$
\begin{array}{l}
\displaystyle t=t(N):= \log{N}/\left\{\widehat{\lambda}_{\partial A}/2+(\lambda_{\Gamma}\vee \lambda_{\Gamma}^{\prime})\right\},
\end{array}
$$
we conclude that
$$
\displaystyle\sup_{s\geq 0}{\EE\left(\Xi_{s}^{\delta/2}\right)^{2/\delta}}\leq c~{N^{-\alpha}}
\quad\mbox{\rm
with}\quad
\alpha=\frac{\widehat{\lambda}_{\partial A}}{\widehat{\lambda}_{\partial A}+2\left(\lambda_{\Gamma}\vee \lambda_{\Gamma}^{\prime}\right)}\in ]0,1].
$$
This ends the proof of the theorem.
\cqfd

\begin{cor}
Assume that $(4^{-1}\sqrt{\lambda_S})\wedge (2^{-1}e\lambda_{R,S})\geq 2$. 
 In this situation,
there exists some $N_0\geq 1$ and some $\alpha\in ]0,1]$ such that for any $N_0\leq N$ and any
initial covariance matrix $P_0$ of the signal  
$$
 \tr(P_0)^2\leq\frac{1}{2}~\frac{\lambda_S}{\lambda_R}
 \left[1+\frac{2}{\lambda_{R}\lambda_S}\right]
\Longrightarrow
\sup_{t\geq 0}{\EE\left(\Vert\xi^1_t-\zeta^1_t\Vert^2\right)}\leq c(P_0)/N^{\alpha}
$$
for some finite constant $ c(P_0)<\infty$ whose value depends on $P_0$.
\end{cor}
\proof
Using (\ref{fv1-3}) we have
$$
\begin{array}{rcl}
d(\xi^1_t-\zeta^1_t)
&=&
\left[(\partial A(m_t)-p_tS)~\xi^1_t+p_tSX_t+A(m_t)-\partial A(m_t)~m_t\right]~dt\\
&&\\
&&-\left[(\partial A(\widehat{X}_t)-P_tS)~\zeta^1_t+P_tSX_t+A(\widehat{X}_t)-\partial A(\widehat{X}_t)~\widehat{X}_t\right]~dt
+d\Ma_t
\end{array}
$$
with the martingale
$$
d\Ma_t:=(p_t-P_t)B^{\prime}R^{-1/2}_{2}d(V_t-\overline{V}^1_{t}).
$$
This yields
$$
\begin{array}{l}
d(\xi^1_t-\zeta^1_t)
\\
\\
=
\left[(\partial A(m_t)-p_tS)~(\xi^1_t-\zeta^1_t)+(p_t-P_t)S(X_t-\zeta^1_t)+(\partial A(m_t)-\partial A(\widehat{X}_t))~\zeta^1_t\right]~dt\\
\\
+\left[\left(A(m_t)-A(\widehat{X}_t)\right)+\left(\partial A(\widehat{X}_t)-\partial A(m_t)\right)~m_t+
\partial A(\widehat{X}_t)~(\widehat{X}_t-m_t)~\right]~dt
+d\Ma_t
\end{array}
$$
with
$$
\sum_{1\leq k\leq r_1}\partial_t\langle\Ma(k),\Ma(k)\rangle_t\leq 2\rho(S)~\Vert~p_t-P_t\Vert^2_F.
$$
This implies that
$$
\begin{array}{l}
d\Vert\xi^1_t-\zeta^1_t\Vert^2\\
\\
\leq 2\langle \xi^1_t-\zeta^1_t,d(\xi^1_t-\zeta^1_t)\rangle+2\rho(S)~\Vert~p_t-P_t\Vert^2_F~dt
\\
\\
\leq \left\{-\lambda_{\partial A}\Vert\xi^1_t-\zeta^1_t\Vert^2+2\rho(S)~\Vert~p_t-P_t\Vert^2_F~+2\Vert\xi^1_t-\zeta^1_t\Vert~\right.\\
\\
\left.\times
\left[\Vert p_t-P_t\Vert_F~\Vert S(X_t-\zeta^1_t)\Vert+\left(\kappa_{\partial A}~\left(\Vert \zeta^1_t\Vert+\Vert m_t\Vert\right)+2\Vert\partial A\Vert\right)~\Vert m_t-\widehat{X}_t\Vert\right]\right\}~dt+d\overline{\Ma}_t
\end{array}
$$
with the martingale
$$
d\overline{\Ma}_t=2\langle \xi^1_t-\zeta^1_t,d\Ma_t\rangle .
$$
Notice that
$$
\begin{array}{l}
\displaystyle2\rho(S)~\Vert~p_t-P_t\Vert^2_F\\
\\
+\displaystyle 2\Vert\xi^1_t-\zeta^1_t\Vert
\left[\Vert p_t-P_t\Vert_F~\Vert S(X_t-\zeta^1_t)\Vert+\left(\kappa_{\partial A}~\left(\Vert \zeta^1_t\Vert+\Vert m_t\Vert\right)+2\Vert\partial A\Vert\right)~\Vert m_t-\widehat{X}_t\Vert\right]\\
\\
\displaystyle\leq \frac{\lambda_{\partial A}}{2}~\Vert\xi^1_t-\zeta^1_t\Vert^2
\times \epsilon_t
\end{array}
$$
with the process
$$
\begin{array}{l}
\epsilon_t
:=\displaystyle2\rho(S)~\Vert~p_t-P_t\Vert^2_F\\
\\
+4\left[\Vert p_t-P_t\Vert_F^2~\Vert S(X_t-\zeta^1_t)\Vert^2+\left(\kappa_{\partial A}~\left(\Vert \zeta^1_t\Vert+\Vert m_t\Vert\right)+2\Vert\partial A\Vert\right)^2~\Vert m_t-\widehat{X}_t\Vert^2\right]/\lambda_{\partial A} .
\end{array}
$$
By Theorem~\ref{theo-c} we have
$$
\sup_{t\geq 0}{\EE(\epsilon_t)}\leq c(P_0)/N^{\alpha}
$$
as soon as $(4^{-1}\sqrt{\lambda_S})\wedge (2^{-1}e\lambda_{R,S})\geq 2$ and
initial covariance matrix $P_0$ of the signal is chosen so that
$$
 \tr(P_0)^2\leq\frac{1}{2}~\frac{\lambda_S}{\lambda_R}
 \left[1+\frac{2}{\lambda_{R}\lambda_S}\right].
$$
This implies that
$$
\partial_t\EE\left(\Vert\xi^1_t-\zeta^1_t\Vert^2\right)
\displaystyle\leq -\frac{\lambda_{\partial A}}{2}~\EE\left(\Vert\xi^1_t-\zeta^1_t\Vert^2\right)+c(P_0)/N^{\alpha} .
$$
The end of the proof of the corollary is now a direct consequence of Gronwall lemma.
\cqfd

\section{Appendix}
\subsection{Regularity conditions}\label{reg-conditions-appendix}
Notice that for any $\alpha,x\geq 0$ we have
 \begin{eqnarray*}
\frac{x}{1+1/x}> 2\alpha&\Longleftrightarrow& x> \alpha\left(1+\sqrt{1+2/\alpha}\right)
 \end{eqnarray*}
 and by \eqref{reference-stability}
 $$
 \lambda_{R,S}
> (8e)^{-1}~{\lambda_{R}~\sqrt{\lambda_{S}}}{~
  \left[1+\frac{1}{\lambda_{R}\sqrt{\lambda_S}}\right]^{-1}}.
 $$
 This shows that
 \begin{eqnarray*}
 (8e)^{-1}~\frac{\lambda_{R}~\sqrt{\lambda_{S}}}{~
  \left[1+\frac{1}{\lambda_{R}\sqrt{\lambda_S}}\right]}>\alpha
  &\Longleftrightarrow& \lambda_{R}~\sqrt{\lambda_{S}}> 4~\alpha~e~\left(1+\sqrt{1+1/(2\alpha e)}\right)\Longrightarrow \lambda_{R,S}> \alpha .
 \end{eqnarray*}  
 Also observe that
 $$
\lambda_{S}>4
\quad\mbox{\rm and}\quad 
 \lambda_{R}~> 2~\alpha~e~\left(1+\sqrt{1+1/(2\alpha e)}\right)\Longrightarrow \lambda_{R,S}> \alpha .
 $$
 This yields the sufficient condition
 $$
  (\lambda_{K}\lambda_{R})\wedge \lambda_{S}>4\quad\mbox{\rm and}\quad
  \lambda_{R}~\sqrt{\lambda_{S}}> 4e~\left(1+\sqrt{1+1/(2e)}\right)
 \Longrightarrow(\ref{reference-stability}) .
 $$  
 Also observe that for any $\alpha\geq 1$ we have
 $$
 \begin{array}{l}
 (\lambda_{K}/\alpha)\wedge (\lambda_{S}/4)
 >1
\quad\mbox{\rm and}\quad 
 \lambda_{R}~> 2~\alpha~e~\left(1+\sqrt{1+1/(2\alpha e)}\right)\\
 \\
 \Longrightarrow
  (\lambda_{K}\lambda_{R}/4)\wedge( \lambda_{R,S}/\alpha)\wedge (\lambda_{S}/4)
 >1.
 \end{array}
 $$
 We end this section with the proof of (\ref{CS-easy-check}). Whenever $\rho(S)=1$ condition (\ref{reference-stability}) 
takes the form
$$
\lambda_{\partial A}>4, \qquad \lambda_{\partial A}^2>4~\kappa_{\partial A}~\tr(R)\quad
\mbox{\rm and}
\quad \lambda_{\partial A}^{3+1/2}>4^2e\left[\tr(R)^2+\frac{1}{2}~\tr(R)~\lambda_{\partial A}^{2}\right]
$$ 
The r.h.s. inequality can be restated as
$$
\left(\frac{\lambda_{\partial A}^{2}}{2}\right)^2\left(1+\frac{1}{4e}~\frac{1}{\sqrt{\lambda_{\partial A}}}\right)>\left(\tr(R)+\frac{\lambda_{\partial A}^{2}}{2}\right)^2
$$
which is equivalent to
$$
\tr(R)<\left(\frac{\lambda_{\partial A}^{2}}{2}\right)\left[\left(1+\frac{1}{4e}~\frac{1}{\sqrt{\lambda_{\partial A}}}\right)^{1/2}-1\right].
$$
This ends the proof of the sufficient condition (\ref{CS-easy-check}).\cqfd
\subsection{Proof of Lemma~\ref{lem-Laplace}}\label{Proof-Laplace-Sec}

We have
\begin{eqnarray*}
-\Gamma_A(t)&=&\lambda_{\partial A}-\left(2\kappa_{\partial A}~\tr(P_t)+\rho(S)~\Vert X_t- \widehat{X}_t\Vert~\right)\leq \lambda_{\partial A}\left[1-2/(\lambda_K\lambda_{R})\right].\end{eqnarray*}
The end of the proof of (\ref{estimation-neg}) is now clear.
Observe that
\begin{eqnarray*}
\Ea_{\Gamma}(t)^{\delta}&=&\exp{\left[\delta\int_0^t\left[
\left(2\kappa_{\partial A}~\tr(P_s)+\rho(S)~\Vert X_s- \widehat{X}_s\Vert~\right)-\lambda_{\partial A}\right]~ ds\right]}\\
&\leq &
\exp{\left[
\delta\lambda_{\partial A}\left[
\frac{2}{\lambda_K}~\left(\tr(P_0)+\frac{1}{\lambda_{R}}\right)-1
\right]~t\right]}~
\exp{\left[\delta~\rho(S)~\int_0^t\Vert X_s- \widehat{X}_s\Vert~ ds\right]}.
\end{eqnarray*}
 We let $\phi_t(x)=X_t$ be the stochastic 
flow of signal starting at $X_0=x$.
We recall the contraction inequality 
\begin{equation}\label{stab-flow}
\Vert \phi_t(x)-\phi_t(y)\Vert\leq \exp{\left(-\lambda_{\partial A}t/2\right)}~\Vert x-y\Vert .
\end{equation}
A proof of (\ref{stab-flow}) can be found in~\cite[Section 3.1]{dkt-1-2016}.
This inequality implies that
\begin{eqnarray*}
\int_0^t~ \Vert X_r-\widehat{X}_r\Vert~dr&=&\int_0^t~ \Vert \phi_{r}(X_0)-\widehat{X}_r\Vert~dr\\
&\leq &\int_0^t~ \Vert \phi_{r}(X_0)-\phi_{r}(\widehat{X}_0)\Vert~dr
+\int_0^t~ \Vert \phi_{r}(\widehat{X}_0)-\widehat{X}_r\Vert~dr\\
&\leq &\left(\int_0^t~e^{-\lambda_{\partial A}r/2}~dr\right)~ \Vert X_0-\widehat{X}_0\Vert+\int_0^t~ \Vert \phi_{r}(\widehat{X}_0)-\widehat{X}_r\Vert~dr\\
&\leq &2\Vert X_0-\widehat{X}_0\Vert/\lambda_{\partial A}+\int_0^t~ \Vert \phi_{r}(\widehat{X}_0)-\widehat{X}_r\Vert~dr.
\end{eqnarray*}
This implies that
$$
\exp{\left[\delta~\rho(S)~\int_0^t\Vert X_s- \widehat{X}_s\Vert~ ds\right]}\leq 
\exp{\left[2\delta\Vert X_0-\widehat{X}_0\Vert/\lambda_S\right]}~
\exp{\left[\delta\rho(S)~\int_0^t~  \Vert \phi_{s}(\widehat{X}_0)-\widehat{X}_s\Vert~ ds\right]}.
$$
Using the estimate $x-1/4\leq x^2$, which is valid for any $x$, we have
$$
\int_0^t~((\Vert \phi_{u}(\widehat{X}_0)-\widehat{X}_u\Vert-1/4)+1/4)~du\leq t/4+\int_0^t~\Vert \phi_{r}(\widehat{X}_0)-\widehat{X}_r\Vert^2dr.
$$
We find that
\begin{eqnarray*}
\exp{\left[\delta\rho(S)~\int_0^t\Vert X_s- \widehat{X}_s\Vert~ ds\right]}&\leq& 
\exp{\left[2 \delta\Vert X_0-\widehat{X}_0\Vert/\lambda_{S}\right]}~\exp{\left(\delta t\rho(S)/4\right)}\\
&&
\times\exp{\left[\delta\rho(S)~\int_0^t~  \Vert \phi_{s}(\widehat{X}_0)-\widehat{X}_s\Vert^2~ ds\right]}.
\end{eqnarray*}
This yields
\begin{eqnarray*}
\EE\left[\exp{\left[\delta~\rho(S)~\int_0^t\Vert X_s- \widehat{X}_s\Vert~ ds\right]}~|~X_0\right]&\leq& ~\exp{\left(t\delta\rho(S)/4\right)}
\exp{\left[2\delta\Vert X_0-\widehat{X}_0\Vert/\lambda_S\right]}\\
&&
\times\EE\left[\exp{\left[\delta\rho(S)~\int_0^t~  \Vert \phi_{s}(\widehat{X}_0)-\widehat{X}_s\Vert^2~ ds\right]}\right).\end{eqnarray*}
We also have the series of inequalities
$$
\begin{array}{l}
\displaystyle\frac{1}{\rho(S)}~\frac{1}{1+\pi_{\partial A}(0)}~\frac{\lambda^2_A}{4\tr(R)}\\
\\
\displaystyle\geq
\frac{1}{\rho(S)}~\frac{\lambda_{\partial A}^2}{4}~\frac{1}{1+\pi_{\partial A}(0)}~\frac{1}{4\tr(R)}
\geq \frac{\lambda_{S}\lambda_R}{2\times 4^2}~\frac{1}{1/2+\tr(P_0)^2(\rho(S)/\tr(R))+\rho(S)\tr(R)/\lambda_{\partial A}^2}~\\
\\
\displaystyle=\frac{1}{4^2}~\lambda_{S}~\lambda_R~\left(2~\frac{\lambda_R}{\lambda_S}~\tr(P_0)^2+\left[1+\frac{2}{\lambda_{R}\lambda_S}\right]\right)^{-1}~\\
\\
\displaystyle\geq \frac{e}{8e}~\sqrt{\lambda_{S}}~\lambda_R~\left[1+\frac{2}{\lambda_{R}\lambda_S}\right]^{-1}~
\left({1+2~\frac{\lambda_R}{\lambda_S}~\tr(P_0)^2\left[1+\frac{2}{\lambda_{R}\lambda_S}\right]^{-1}}\right)^{-1}\\
\\
\displaystyle=e~\lambda_{R,S}~\left(1+2~\frac{\lambda_R}{\lambda_S}~\tr(P_0)^2\left[1+\frac{2}{\lambda_{R}\lambda_S}\right]^{-1}\right)^{-1}.
\end{array}
$$
This shows that
$$
\delta\rho(S)\leq 
\frac{\epsilon}{1+\pi_{\partial A}(0)}~\frac{\lambda^2_A}{4\tr(R)}
$$
for some $\epsilon\in [0,1]$ as soon as
$$
\tr(P_0)^2\leq \frac{1}{2}~\frac{\lambda_S}{\lambda_R}
 \left[1+\frac{2}{\lambda_{R}\lambda_S}\right]\left(
\frac{e}{\delta}~\epsilon~\lambda_{R,S}-1\right)
\quad\mbox{\rm
for any $\delta\leq e~\epsilon~\lambda_{R,S}$} .
$$
The end of the proof of (\ref{estimation-positive}) is a direct consequence of~\cite[Theorem 3.2]{dkt-1-2016}.

The last assertion resumes to~\cite[Lemma 4.1]{dkt-1-2016}.
This ends the proof of the lemma.
\cqfd

\subsection{Proof of Lemma~\ref{prop-unif-control}}\label{Lemma-Unif-control-sec}

Using (\ref{f23}) we have
\begin{eqnarray*}
d\tr(p_t)&=&\left(\tr((\partial A\left[m_t\right]+\partial A\left[m_t\right]^{\prime})p_t)-\tr(Sp^2_t)+\tr(R)\right)~dt+\frac{1}{\sqrt{N-1}}~d\Ma_t\\
&\leq &\left[-\lambda_{\partial A}~\tr(p_t)-r_1^{-1}\rho(S)~\tr(p_t)^2+\tr(R)\right]~dt+\frac{1}{\sqrt{N-1}}~d\Ma_t\
\end{eqnarray*}
with a martingale $\Ma_t$ with an angle bracket
$$
\partial_t\langle\Ma\rangle_t=4\tr((R+p_tSp_t)p_t)\leq 4\tr(p_t)~\left(\rho(R)+\rho(S)~\tr(p_t)^2\right).
$$
Using~\cite[Lemma 4.1]{dt-1-2016} we have
$$
1\leq n\leq 1+\frac{(N-1)}{2r_1}\frac{\rho(S)}{\lambda_{\tiny max}(S)}\Longrightarrow
\sup_{t\geq 0}{\EE\left(\tr(p_t)^n\right)}<\infty .
$$
By (\ref{f21}) we have
$$
dm_t=[A\left[m_t\right]-p_tSm_t+p_t SX_t ]~dt+p_t~B^{\prime}R^{-1/2}_{2}~dV_t+\frac{1}{\sqrt{N}}~d\overline{M}_t.
$$
Since $\overline{M}_t$ is independent of $V_t$ we have
$$
d\Vert m_t\Vert^2=\left(2~\langle m_t,[A\left[m_t\right]-p_tSm_t+p_t SX_t ]\rangle+\tr(R+p_tSp_t)\right)~dt+
d\widetilde{M}_t
$$
with the martingale
$$
d\widetilde{M}_t=2~\langle m_t,p_t~B^{\prime}R^{-1/2}_{2}~dV_t\rangle+2~\frac{1}{\sqrt{N}}~\langle m_t,~d\overline{M}_t\rangle
$$
and the angle bracket
\begin{eqnarray*}
\partial_t\langle \widetilde{M}\rangle_t&=&4~\langle m_t, (R+p_tSp_t)m_t\rangle/N+
4~\langle m_t,(p_tSp_t)m_t\rangle\leq \Va_t~\Vert m_t\Vert^2
\end{eqnarray*}
with
$$
\Va_t~:=4~[\tr(R+p_tSp_t)/N+\tr(p_tSp_t)].
$$
Observe that 
\begin{eqnarray*}
\langle m_t,A\left[m_t\right]\rangle&=&\langle m_t-0,A\left[m_t\right]-A(0)\rangle+\langle m_t,A\left[0\right]\rangle\\
&\leq &-\lambda_{A}~\Vert m_t\Vert^2+\Vert A(0)\Vert~\Vert m_t\Vert\leq -(\lambda_{A}/2)~\Vert m_t\Vert^2+\Vert A(0)\Vert^2/(2\lambda_A).
\end{eqnarray*}
This yields the estimate
\begin{eqnarray*}
d\Vert m_t\Vert^2&\leq& \left(-\lambda_{A}~\Vert m_t\Vert^2+\Vert A(0)\Vert^2/\lambda_A+2\Vert m_t\Vert~\Vert p_tS\Vert~\Vert X_t\Vert
+\tr(R+p_tSp_t)\right)~dt+
d\widetilde{M}_t
\end{eqnarray*}
from which we find that
$$
d\Vert m_t\Vert^2\leq  \left(-\frac{\lambda_{A}}{2}~\Vert m_t\Vert^2+\Ua_t\right)~dt+\sqrt{\Va_t}~
d\Na_t
$$
with $\partial_t\langle \Na\rangle_t\leq 1$ and
\begin{eqnarray*}
\Ua_t&:=&\Vert A(0)\Vert^2/\lambda_A+\Vert p_tS\Vert^2~\Vert X_t\Vert^2/\lambda_A
+\tr(R+p_tSp_t).
\end{eqnarray*}
Arguing as in the proof of Theorem~\ref{theo-stability} we conclude that
$$
\forall 1\leq 3n\leq 1+{(N-1)}/{(2r_1)}\qquad \sup_{t\geq 0}{\EE\left(\Vert m_t\Vert^{2n}\right)}<\infty .
$$
Using (\ref{fv1-3}) we have
$$
d\xi^1_t=\left[(\partial A\left[m_t\right]-p_tS)~\xi^1_t+p_tSX_t+A\left[m_t\right]-\partial A\left[m_t\right]~m_t\right]~dt+d\Ma_t
$$
with the martingale
$$
d\Ma_t:=R^{1/2}_{1}d\overline{W}_t^1+p_tB^{\prime}R^{-1/2}_{2}d(V_t-\overline{V}^1_{t}).
$$
This implies that
$$
\begin{array}{rcl}
d\Vert \xi^1_t\Vert^2&=&\left[2\langle\xi^1_t,
\left[(\partial A\left[m_t\right]-p_tS)~\xi^1_t+p_tSX_t+A\left[m_t\right]-\partial A\left[m_t\right]~m_t\right]\rangle\right.\\
\\
&&\hskip6cm\left.+\tr(R)+2\tr(p_tSp_t)\right]~dt+d\overline{\Ma}_t\\
\\
&\leq &\left[-(\lambda_{\partial A}/2)~\Vert \xi^1_t\Vert^2+\Ua_t\right]~dt+d\overline{\Ma}_t
\end{array}
$$
with
$$
d\overline{\Ma}_t:=2\langle \xi^1_t,d\Ma_t\rangle\Longrightarrow\partial_t\langle \overline{\Ma}\rangle_t\leq \Va_t~\Vert\xi^1_t\Vert^2~
$$
and
\begin{eqnarray*}
\Ua_t&:=&2\Vert p_tSX_t+A\left[m_t\right]-\partial A\left[m_t\right]~m_t\Vert^2/\lambda_{\partial A}+\tr(R)+2\tr(Sp^2_t),\\
\Va_t&=&4\left(\tr(R)+2\tr(p_tSp_t)\right).
\end{eqnarray*}
The end of the proof follows the same arguments as above, so it is skipped.
This completes the proof of the lemma.
\cqfd

\subsection{Proof of Lemma~\ref{lem-perturbation}}\label{Proof-perturbation-lemma}

By (\ref{f21}) and   \eqref{f23} we have
$$
d(p_t-P_t)=\Pi_t~dt+d\Ma_t\quad\mbox{\rm and}\quad
d(m_t-\widehat{X}_t)
=\overline{\Pi}_t~dt+d\overline{\Ma}_t
$$
with the drift terms
$$
\begin{array}{rcl}
\Pi_t&=&
\left(
\partial A(m_t)p_t-\partial A(\widehat{X}_t)P_t
\right)+\left(
\partial A(m_t)p_t-\partial A(\widehat{X}_t)P_t
\right)^{\prime}\\
&&\\
&&\hskip3cm+(P_t-p_t)SP_t+((P_t-p_t)SP_t)^{\prime}-(P_t-p_t)S(P_t-p_t)\\
&&\\
\overline{\Pi}_t&=&(A(m_t)-A(\widehat{X}_t))-p_tS(m_t-\widehat{X}_t)+(p_t-P_t) S(X_t- \widehat{X}_t)
\end{array}
$$
and the martingales
$$
d\Ma_t:=\frac{1}{\sqrt{N-1}}~dM_t , \qquad
d\overline{\Ma}_t:=(p_t-P_t)~B^{\prime}R^{-1/2}_{2}~dV_t+
\frac{1}{\sqrt{N}}~d\overline{M}_t.
$$
Using the decomposition
$$
\partial A(m_t)p_t-\partial A(\widehat{X}_t)P_t=\partial A(m_t)(p_t-P_t)+(\partial A(m_t)-\partial A(\widehat{X}_t))P_t,
$$
we check that
$$
\begin{array}{l}
\Pi_t
=\left[\partial A(m_t)-\frac{1}{2}(p_t+P_t)S\right](p_t-P_t)+(p_t-P_t)\left[\partial A(m_t)-\frac{1}{2}(p_t+P_t)S\right]^{\prime}\\
\\
\hskip4cm+(\partial A(m_t)-\partial A(\widehat{X}_t))P_t+P_t(\partial A(m_t)-\partial A(\widehat{X}_t))^{\prime}.
\end{array}
$$
This implies that
\begin{eqnarray*}
\langle p_t-P_t,\Pi_t\rangle&\leq& -\lambda_{\partial A}~\Vert p_t-P_t\Vert_F^2+2\kappa_{\partial A}~\tr(P_t)~\Vert p_t-P_t\Vert_F~\Vert m_t-\widehat{X}_t\Vert
\end{eqnarray*}
from which we prove that
\begin{eqnarray*}
d\Vert p_t-P_t\Vert_F^2&=&2~\langle p_t-P_t,d(p_t-P_t)\rangle\\
&&+\frac{2}{N-1}~\left[\tr((R+p_tSp_t)p_t)+\tr(R+p_tSp_t)\tr(p_t)\right]~dt\\
&\leq &\left\{-2\lambda_{\partial A}~\Vert p_t-P_t\Vert_F^2+4\kappa_{\partial A}~\tr(P_t)~\Vert p_t-P_t\Vert_F~\Vert m_t-\widehat{X}_t\Vert\right.\\
&&\hskip.3cm\left.+\frac{2}{N-1}~\left[\tr((R+p_tSp_t)p_t)+\tr(R+p_tSp_t)\tr(p_t)\right]\right\}~dt+d\Na_t
\end{eqnarray*}
with the martingale
$$
d\Na_t=\frac{2}{\sqrt{N-1}}~\langle p_t-P_t,dM_t\rangle=\tr((p_t-P_t)dM_t).
$$
After some computations we find that
$$
\partial_t\langle \Na\rangle_t\leq \frac{4}{N-1}~\Vert p_t-P_t\Vert_F^2~\tr(p_t(R+p_tSp_t)).
$$
In much the same vein we have
\begin{eqnarray*}
\langle m_t-\widehat{X}_t,\overline{\Pi}_t\rangle&=&\langle m_t-\widehat{X}_t,(A(m_t)-A(\widehat{X}_t))-p_tS(m_t-\widehat{X}_t)+(p_t-P_t) S(X_t- \widehat{X}_t)\rangle\\
&\leq &-\lambda_A~\Vert m_t-\widehat{X}_t\Vert^2+\rho(S)~\Vert X_t- \widehat{X}_t\Vert~\Vert p_t-P_t\Vert_F~\Vert m_t-\widehat{X}_t\Vert .
\end{eqnarray*}
This implies that
\begin{eqnarray*}
d\Vert m_t-\widehat{X}_t\Vert^2&=&2~\langle (m_t-\widehat{X}_t),d(m_t-\widehat{X}_t)\rangle\\
&&+\left(\tr(S(p_t-P_t)^2)+\frac{1}{N}~\tr(R+p_tSp_t)\right)~dt\\
&\leq &\left\{-2\lambda_A~\Vert m_t-\widehat{X}_t\Vert^2+2\rho(S)~\Vert X_t- \widehat{X}_t\Vert~\Vert p_t-P_t\Vert_F~\Vert m_t-\widehat{X}_t\Vert\right.\\
&&\left. +\rho(S)\Vert p_t-P_t\Vert_F^2+\frac{1}{N}~\tr(R+p_tSp_t)\right\} dt +d\overline{\Na}_t
\end{eqnarray*}
with the martingale
\begin{eqnarray*}
d\overline{\Na}_t=2~\langle (m_t-\widehat{X}_t),d\overline{\Ma}_t\rangle=
2~\langle (m_t-\widehat{X}_t),(p_t-P_t)~B^{\prime}R^{-1/2}_{2}~dV_t\rangle+\frac{2}{\sqrt{N}}~\langle (m_t-\widehat{X}_t),d\overline{M}_t\rangle .
\end{eqnarray*}
In addition we have
\begin{eqnarray*}
\partial_t\langle \overline{\Na}\rangle_t&\leq& 4\rho(S)~\Vert m_t-\widehat{X}_t\Vert^2\Vert p_t-P_t\Vert_F^2+\frac{4}{N}~\langle (m_t-\widehat{X}_t),(R+p_tSp_t)(m_t-\widehat{X}_t)\rangle\\
&\leq &2\rho(S)~\left(\Vert m_t-\widehat{X}_t\Vert^2+\Vert p_t-P_t\Vert_F^2\right)^2+\frac{4}{N}~\Vert m_t-\widehat{X}_t\Vert^2~\tr(R+p_tSp_t).
\end{eqnarray*}

Combining the above estimates we find that
$$
\begin{array}{l}
d\,\Xi_t
\leq \left\{-2\lambda_A~\Vert m_t-\widehat{X}_t\Vert^2+2~\Vert p_t-P_t\Vert_F~\Vert m_t-\widehat{X}_t\Vert~\left(2\kappa_{\partial A}~\tr(P_t)+\rho(S)~\Vert X_t- \widehat{X}_t\Vert\right)\right\}~dt\\
\\
\hskip3cm\displaystyle-\left(2\lambda_{\partial A}-\rho(S)\right)\Vert p_t-P_t\Vert_F^2~dt\\
\\
\hskip.3cm\displaystyle+\frac{1}{N}\left\{\tr(R+p_tSp_t)+\frac{2N}{N-1}~\left[\tr((R+p_tSp_t)p_t)+\tr(R+p_tSp_t)\tr(p_t)\right]~\right\}~dt+d\Na_t+d\overline{\Na}_t.\ 
\end{array}
$$
Recalling that  $$
2\lambda_A\geq \lambda_{\partial A}> 0\quad\mbox{\rm and}\quad
2\lambda_{\partial A}- \rho(S)\geq \lambda_{\partial A},
$$ 
this yields the estimate
\begin{eqnarray*}
d\,\Xi_t
&\leq &\left\{-\lambda_{\partial A}~\Xi_t+2~\Vert p_t-P_t\Vert_F~\Vert m_t-\widehat{X}_t\Vert\left(2\kappa_{\partial A}~\tr(P_t)+\rho(S)~\Vert X_t- \widehat{X}_t\Vert~\right)\right\}~dt\\
&&+\frac{1}{N}\left(\frac{4N}{N-1}~\tr(R+p_tSp_t)~\tr(p_t)+~\tr(R+p_tSp_t)\right)~dt+d\overline{\Na}_t+d\Na_t.
\end{eqnarray*}
On the other hand using the inequality $2ab\leq a^2+b^2$ we prove that
\begin{eqnarray*}
d\,\Xi_t
&\leq &-\left\{\lambda_{\partial A}-\left(2\kappa_{\partial A}~\tr(P_t)+\rho(S)~\Vert X_t- \widehat{X}_t\Vert~\right)\right\}~\Xi_t~dt\\
&&+\frac{1}{N}~\left(1+\frac{4}{1-1/N}~\tr(p_t)\right)~\left[\tr(R)+\tr(S)\tr(p_t)^2\right]~dt+d\overline{\Na}_t+d\Na_t
\end{eqnarray*}
from which we conclude that 
\begin{eqnarray*}
d\,\Xi_t
&\leq &\left(\Gamma_A(t)~\Xi_t
~+\frac{1}{N}~\Ua_t\right)dt+d\,\Upsilon_t\quad\mbox{\rm with}\quad \Ua_t:=(1+8\tr(p_t))~\left[\tr(R)+\rho(S)~\tr(p_t)^2\right]
\end{eqnarray*}
 and the martingale $\Upsilon_t:=\Upsilon_t^{(1)}+\Upsilon_t^{(2)}$ given by
\begin{eqnarray*}
d\,\Upsilon_t^{(1)}&:=&
2~\langle (m_t-\widehat{X}_t),(p_t-P_t)~B^{\prime}R^{-1/2}_{2}~dV_t\rangle\\
d\,\Upsilon_t^{(2)}&:=&\frac{2}{\sqrt{N}}~\langle (m_t-\widehat{X}_t),d\overline{M}_t\rangle+\frac{2}{\sqrt{N-1}}~\langle p_t-P_t,dM_t\rangle .
\end{eqnarray*}
Observe that
\begin{eqnarray*}
\langle \Upsilon^{(1)},\Upsilon^{(2)}\rangle_t&=&0\\
\partial_t\langle \Upsilon^{(1)}\rangle_t&\leq& 
2\rho(S)~\left(\Vert m_t-\widehat{X}_t\Vert^2+\Vert p_t-P_t\Vert_F^2\right)^2\leq 2\rho(S)~\Xi_t^2\\
\partial_t\langle \Upsilon^{(2)}\rangle_t&\leq& \frac{4}{N}\left[~\Vert m_t-\widehat{X}_t\Vert^2~+2~\Vert p_t-P_t\Vert_F^2~\tr(p_t)\right]~\tr(R+p_tSp_t)\leq  \frac{4}{N}~\Ua_t~\Xi_t .\end{eqnarray*}

This ends the proof of the lemma.\cqfd

\subsubsection*{Acknowledgements}
We would like to thank an anonymous referee for the solving of a technical issue in the earlier proof of Theorem 3.2. Her/His detailed comments also greatly improved the 
presentation of the article.

Pierre DEL MORAL\\
INRIA Bordeaux Research Center (France) \& UNSW School of Mathematics and Statistics  (Australia)\\
p.del-moral@unsw.edu.au
\medskip

Aline KURTZMANN\\
Universit\'e de Lorraine, Institut Elie Cartan de Lorraine \\ CNRS, Institut Elie Cartan de Lorraine, UMR 7502, Vandoeuvre-l\`es-Nancy, F-54506, France.\\
aline.kurtzmann@univ-lorraine.fr
\medskip

Julian TUGAUT\\
Univ Lyon, Universit\'e Jean Monnet, CNRS UMR 5208, Institut Camille Jordan\\
 Maison de l'Universit\'e, 10 rue Tr\'efilerie, CS 82301, 42023 Saint-Etienne Cedex 2, France\\
 tugaut@math.cnrs.fr
 
\end{document}